\documentclass{article}

\usepackage{PRIMEarxiv}

\usepackage[utf8]{inputenc} 
\usepackage[T1]{fontenc}    
\usepackage{hyperref}       
\usepackage{url}            
\usepackage{booktabs}       
\usepackage{amsfonts}       
\usepackage{nicefrac}       
\usepackage{microtype}      
\usepackage{lipsum}
\usepackage{fancyhdr}       
\usepackage{graphicx}       
\graphicspath{{media/}}     

\usepackage{amsmath}
\usepackage{amsthm}
\usepackage{enumerate}
\usepackage{comment}

\title{On the General Position Number of Cacti
}

\author{
  Toshihiko Kimura \\
  Graduate School of Engineering \\
  Nagoya Institute of Technology\\
  Aichi, Japan\\
  \texttt{t.kimura.185@stn.nitech.ac.jp} \\
   \And
  Yonghwan Kim \\
  Graduate School of Engineering \\
  Nagoya Institute of Technology\\
  Aichi, Japan\\
  \texttt{kim@nitech.ac.jp} \\
   \And
  Yoshiaki Katayama \\
  Graduate School of Engineering \\
  Nagoyda Institute of Technology\\
  Aichi, Japan\\
  \texttt{katayama@nitech.ac.jp} \\
}

\newcommand{\GPS}{\mathtt{GPS}}
\newcommand{\GPP}{\mathtt{GPP}}
\newcommand{\GPN}{\mathtt{GPN}}

\newcommand{\sGPS}{\mathsf{s}\text{-}\mathtt{GPS}}
\newcommand{\subGPS}{\mathsf{sub}\text{-}\mathtt{GPS}}

\newcommand{\defeq}{\overset{\text{\tiny def}}{=}}

\newcommand{\ie}{{\it i.e.,~}}

\newcommand{\etal}{{\it et al.}}

\newcommand{\bd}{\emph{beard}} 
\newcommand{\bc}{\emph{bearded cycle}} 
\newcommand{\nt}{\emph{non-beard tree}}

\newtheorem{theorem}{Theorem}[section]
\newtheorem{lemma}[theorem]{Lemma}
\theoremstyle{definition}
\newtheorem{definition}[theorem]{Definition}
\newtheorem*{theorem*}{Theorem}

\begin{document}
\maketitle

\begin{abstract}
Given a graph $G=(V,E)$, where $V$ is the set of vertices and $E$ is the set of edges, 
a subset of vertices $S \subseteq V$ is called a {\it General Position Set} ($\GPS$) 
if no shortest path on $G$ between any two vertices in $S$ contains another vertex from $S$.

The {\it General Position Number} ($\GPN$) of a graph $G$ 
is the maximum size of a $\GPS$ in $G$. 
The decision version of this problem, 
determining whether the $\GPN$ of a graph $G$ is at least a given integer $i$,
is known as the {\it General Position Problem} ($\GPP$). 
The $\GPP$ is known to be \emph{NP-complete}, 
making the computation of the $\GPN$ for general graphs computationally challenging. 
Hence, much of the existing research has focused on determining the $\GPN$ 
for special classes of graphs. 
One such class is that of cactus graphs (\ie cacti), 
in which any two simple cycles share at most one vertex and no edge.
For cactus graphs, several upper and lower bounds on the $\GPN$ have been previously presented. 

In this paper, we determine the exact value of the $\GPN$ for cactus graphs. 
We propose a novel graph decomposition that partitions a cactus graph 
into two types of subgraphs, called \bc s and \nt s.
For each subgraph, we compute a special $\GPS$, named a $\subGPS$.
We prove that the union of all maximum $\subGPS$s constitutes a maximum $\GPS$ 
for the original cactus graph.
This result allows us to compute the exact $\GPN$ for any cactus graph.
\end{abstract}



\section{Background}
\label{sec:background}
In graph theory, the concept of a {\it General Position Set} ($\GPS$) was firstly introduced 
by Manuel \etal~in 2018 \cite{MANUEL_KLAVŽAR_2018}. 
Given a graph $G=(V,E)$, a subset $S \subseteq V$ is called a $\GPS$ 
if no shortest path between any two vertices in $S$ contains another vertex from $S$.

$\GPS$s have potential applications in efficient and secure communication in networks. 
In a distributed system consisting of multiple computational entities, 
the system can be modeled as a graph, 
where vertices represent computing terminals and edges represent communication links between them. 
If the terminals form a $\GPS$, then any shortest path between them avoids all
other terminals in the set, enabling secure point-to-point communication. 
In such networks, it is desirable for the secure set of terminals to be as large as possible.

A fundamental problem related to $\GPS$s is the {\it General Position Problem} ($\GPP$). 
Given a graph $G$ and an integer $i$, the $\GPP$ asks whether the size of the largest $\GPS$ in $G$,
known as the {\it General Position Number} ($\GPN$), is at least $i$. 
It has been shown that the $\GPP$ is \emph{NP-complete} for general graphs \cite{MANUEL_KLAVŽAR_2018},
making the computation of the $\GPN$ computationally intractable in the general case.
Therefore, research has focused on analyzing restricted classes of graphs 
\cite{anand2019characterizationgeneralpositionsets,klavžar2021generalpositionsetscartesian,klavžar2020generalpositionnumberinteger,yao2020generalpositionsetclasses}. 

For cactus graphs, the focus of this study, 
previous work has provided upper and lower bounds on the $\GPN$ for general cacti, 
as well as exact values for certain subclasses \cite{yao2020generalpositionsetclasses}. 
This study aims to determine the exact value of the $\GPN$ for all cactus graphs.


\section{Related Work}

The \emph{no-three-in-line problem} is a classical combinatorial problem 
that seeks the largest subset of vertices on a two-dimensional grid 
such that no three points lie on a straight line. 
One of the earliest instances of this problem was introduced by Dudeney in the puzzle 
\emph{A Puzzle with Pawns} \cite{dudeney1958amusements}, which asks how to place 14 additional 
emph{pawns} on a chessboard,
given that two are already placed, such that no three pawns are aligned. 
Variants of this problem have been studied on toroidal grids and in three-dimensional grids \cite{fowler2012nothreeinlineproblemtorus,10.1007/978-3-540-31843-9_40}.

Inspired by this problem, 
Froese \etal~ introduced the \emph{General Position Subset Selection Problem}, 
which asks whether, given a set $P$ of points in the plane and a positive integer $k$, 
there exists a subset $S \subseteq P$ of size at least $k$ such that no three points in $S$ 
are collinear \cite{froese2017findingpointsgeneralposition}. 
They showed that this problem is \emph{NP-hard} \cite{froese2017findingpointsgeneralposition}.

Building on these foundational problems, Manuel \etal~ formulated the {\it General Position Problem} ($\GPP$) 
in the context of graphs \cite{MANUEL_KLAVŽAR_2018}. 
A {\it General Position Set} ($\GPS$) in a graph $G=(V,E)$ is defined as a subset $S \subseteq V$ 
such that no shortest path between any two vertices in $S$ contains another vertex from $S$ \cite{MANUEL_KLAVŽAR_2018}. 
The maximum size of a $\GPS$ is called the {\it General Position Number} ($\GPN$),
and the $\GPP$ asks whether the $\GPN$ of a given graph $G$ is at least a given integer $i$ \cite{MANUEL_KLAVŽAR_2018}. 
This problem has been shown to be \emph{NP-complete} for general graphs, 
making the computation of the $\GPN$ computationally challenging in most cases. 
As a result, research has shifted toward studying restricted graph classes \cite{anand2019characterizationgeneralpositionsets,klavžar2021generalpositionsetscartesian,klavžar2020generalpositionnumberinteger,yao2020generalpositionsetclasses}. 

For cacti, the focus of this study, 
prior work has established upper and lower bounds on the $\GPN$ \cite{yao2020generalpositionsetclasses}. 
Klavžar \etal~ determined the exact $\GPN$ for two-dimensional grid graphs under certain conditions, 
enumerated the number of maximum $\GPS$ in arbitrary two-dimensional grids, 
provided bounds for three-dimensional grids, and derived exact values for $n$-dimensional grids \cite{klavžar2021generalpositionsetscartesian,klavžar2020generalpositionnumberinteger}. 
Anand \etal~ also computed the exact $\GPN$ for the complements of certain bipartite graphs \cite{anand2019characterizationgeneralpositionsets}.

\begin{figure}[tb]
    \begin{minipage}{0.49\hsize}
        \begin{center}
            \includegraphics[width=0.5\linewidth]{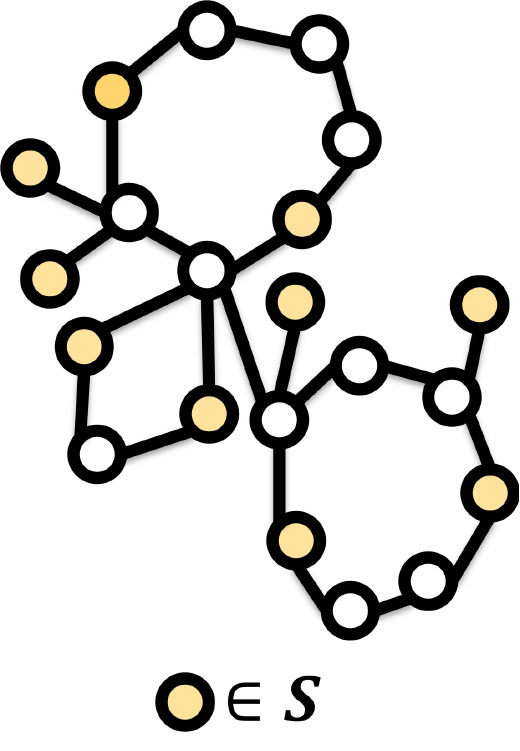}
        \end{center}
        \caption{Case that all cycles are good and $\subGPS$ $S$}
        \label{fig:goodcycles}
    \end{minipage}
\hfill
    \begin{minipage}{0.49\hsize}
        \begin{center}
            \includegraphics[width=0.5\linewidth]{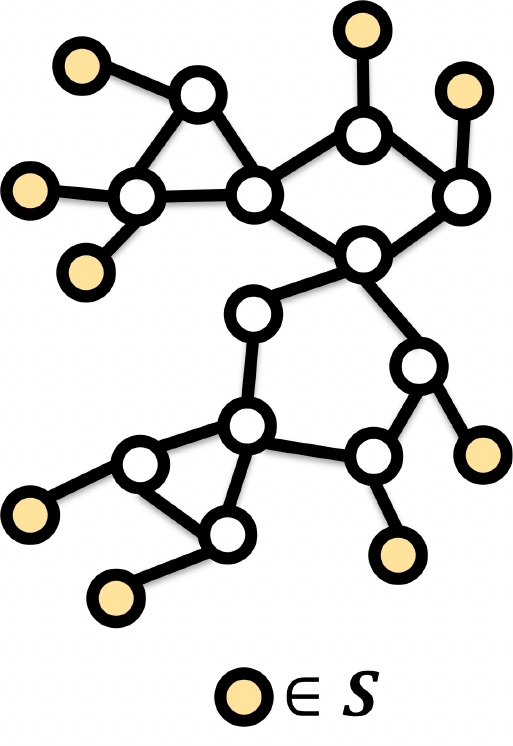}
        \end{center}
        \caption{Case that all cycles are bad and $\subGPS$ $S$}
        \label{fig:badcycles}
        \end{minipage}
\end{figure}

Yan \etal~ proposed upper and lower bounds on the $\GPN$ for general cactus graphs 
and provided exact values for certain subclasses \cite{yao2020generalpositionsetclasses}. 
Let $C_n^{t,k}$ denote the class of cactus graphs with $n$ vertices, $k$ simple cycles, and $t$ leaves. 
For any graph $G \in C_n^{t,k}$, let $C_l=(u_1,u_2,\dots,u_l,u_1)=(V_{C_l},E_{C_l})$ be a simple cycle in $G$. 
For any pair $u_i, u_j \in V_{C_l}$, the cycle $C_l$ contains two simple paths connecting them. 
Among these, let $d_C(u_i,u_j)$ denote the number of edges in the path that includes all articulation points of $G$ lying on $C_l$. 
Define $D_C = \min_{u_i,u_j \in V_{C_l}} d_C(u_i,u_j)$. A cycle $C_l$ is considered \emph{good} if either 
(i) $l$ is odd and $D_C \le \lfloor \frac{l}{2} \rfloor - 1$, 
or (ii) $l$ is even and $D_C \le \frac{l}{2} - 2$. 

The upper bound on the $\GPN$ of $G$ is given by $\max(\{3, 2k + t\})$, with equality if and only if 
all simple cycles in $G$ are good. The lower bound is $t$, with equality if and only if none of the simple cycles in $G$ 
are good (see Figures \rm{\ref{fig:goodcycles}} and \rm{\ref{fig:badcycles}}). 

As highlighted above, considerable research has focused on determining the $\GPN$ for specific graph classes.
However, the exact value of the $\GPN$ for general cactus graphs remains unknown.
This study seeks to address this gap by deriving the exact $\GPN$ for general cactus graphs.


\section{Preliminaries}
In this section, we define several subgraphs and vertex sets of a cactus graph
studied in this paper.

We consider simple, undirected graphs $G = (V, E)$,
where $V$ is the set of vertices, and $E \subseteq V \times V$ is the set of edges. 
For two different vertices $v_i, v_j \in V$ ($v_i \ne v_j$), if $(v_i, v_j) \in E$, 
we say that $v_i$ and $v_j$ are adjacent. 
We assume there are no self-loops; that is, $(v_i, v_j) \notin E$ for $v_i = v_j$.

We begin with the definition of a cactus graph, the class of graphs on which this study focuses.

\begin{definition} \rm{(Cactus Graph)}
A connected graph $G$ is called a \emph{cactus graph} 
if no two simple cycles in $G$ share an edge.
\end{definition}

\begin{figure}[tb]
    \begin{center}
        \begin{minipage}[htpb]{0.45\columnwidth}
        \begin{center}
            \includegraphics[width=0.5\linewidth]{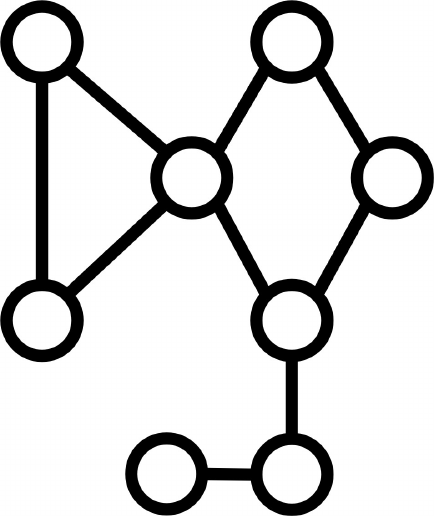}
        \end{center}
            \caption{A example of a cactus graph}\label{fig:exp_yescactus}
        \end{minipage}   
        \hfill
        \begin{minipage}[htbp]{0.54\columnwidth}
        \begin{center}
            \includegraphics[width=0.43\linewidth]{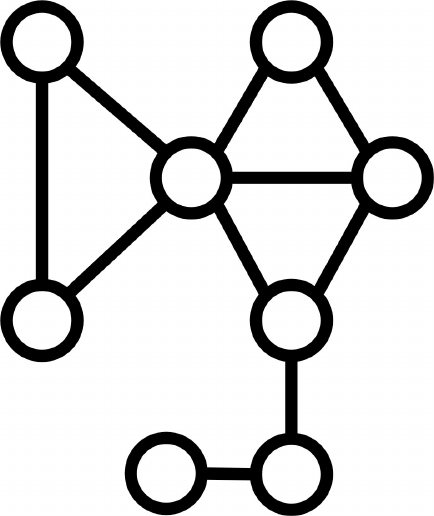}
        \end{center}
            \caption{An example of a non-cactus graph}\label{fig:exp_nocactus}
        \end{minipage}
    \end{center}
\end{figure}

Figures \rm{\ref{fig:exp_yescactus}} and \rm{\ref{fig:exp_nocactus}} illustrate examples of a cactus graph 
and a non-cactus graph, respectively.

Next, we define the concepts of $\bd$, $\bc$, and $\nt$.
These components are used to decompose a cactus graph in order to determining its $\GPN$.

\begin{definition} \rm{(Beard and Non-beard Tree)}
Let $G$ be a cactus graph, and let $G_F$ be the graph obtained by removing all edges included in cycles 
and all degree-2 vertices belonging to cycles. 
A connected component of $G_F$ that contains exactly two leaves on $G_F$, 
one of which belongs to exactly one cycle and the other not belonging to any cycle, is called a $\bd$.
A connected component of $G_F$ that is not a beard is called a $\nt$.
\end{definition}

\begin{figure}[tb]
    \begin{center}
    \includegraphics[width=0.6\linewidth]{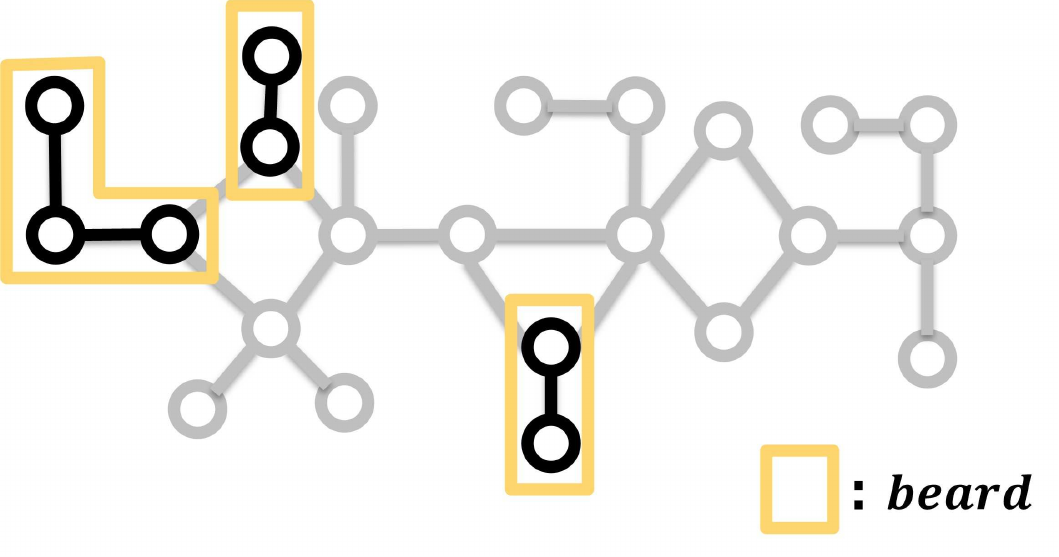}
    \caption{An example of beards}\label{fig:exp_beard}
    \end{center}
\end{figure}

\begin{figure}[tb]
    \begin{center}
    \includegraphics[width=0.6\linewidth]{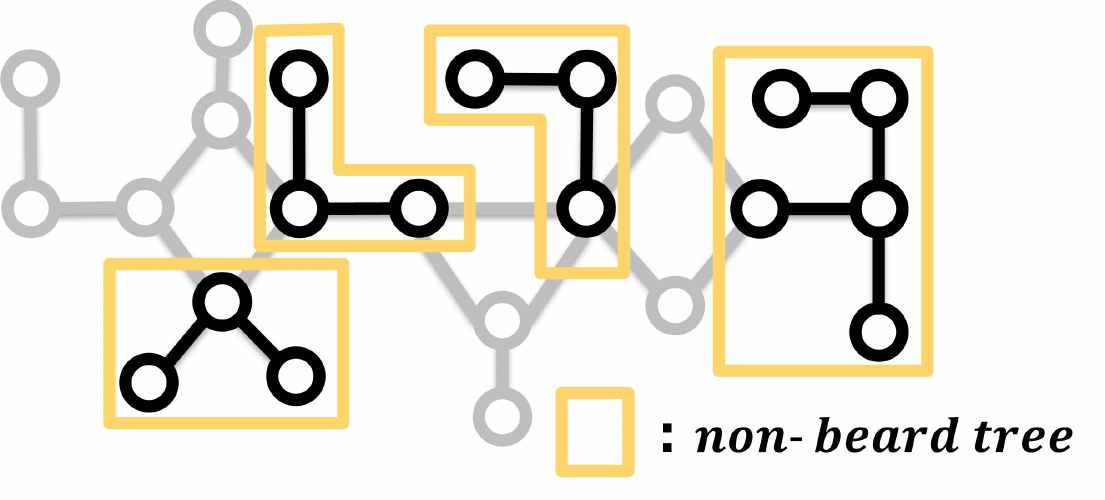}
    \caption{An example of non-beard trees}\label{fig:exp_nonbeardtree}
    \end{center}
\end{figure}

Examples of beards and \nt s are shown in Figures \rm{\ref{fig:exp_beard}} and \rm{\ref{fig:exp_nonbeardtree}}, respectively.

\begin{figure}[tb]
    \begin{center}
    \includegraphics[width=0.6\linewidth]{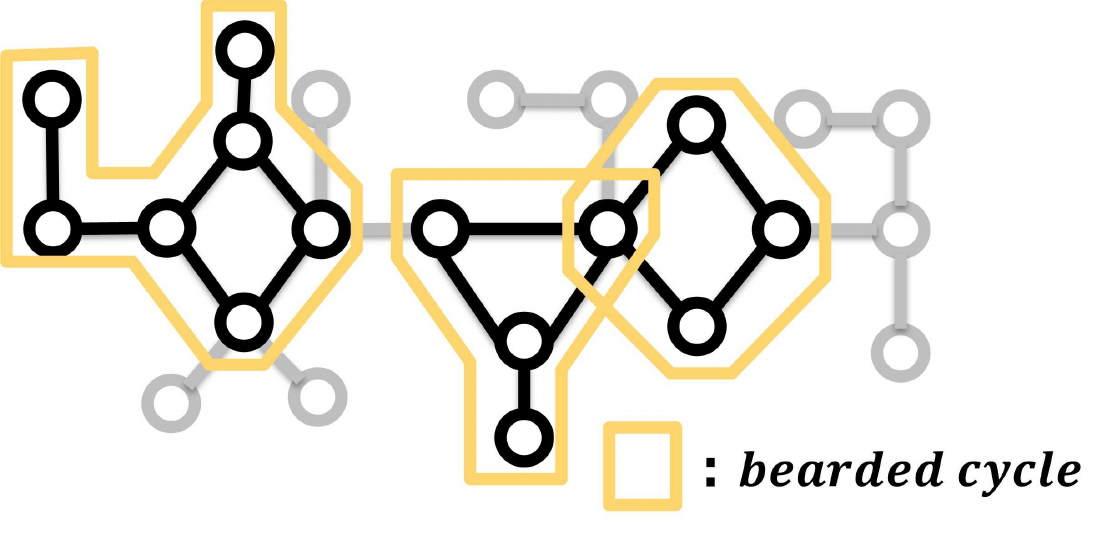}
    \caption{Example of bearded cycles}\label{fig:exp_beardcycle}
    \end{center}
\end{figure}

\begin{definition} \rm{(Bearded Cycle)}
For each cycle $G_C$ in a cactus graph, 
the graph consisting of $G_C$ and all $\bd$s sharing a vertex with $G_C$ is called a $\bc$.
\end{definition}

An example of a bearded cycle is shown in Figure \rm{\ref{fig:exp_beardcycle}}.

We now define several important vertex sets used in this paper. 
Let:
\begin{itemize}
\item $L_G$ be the set of all leaves in $G$.
\item $A_G$ be the set of all articulation points in $G$.
\item $C_G$ be the set of all vertices contained in any cycle in $G$.
\item $B_G$ be the set of all vertices contained in any beard in $G$. 
\end{itemize}

For any vertex subset $S \subseteq V$, 
let $V_{mP}(G,S)$ denote the set of all vertices lying on any shortest paths between any pair of vertices $u,v \in S$ in $G$.

For any bearded cycle $G_{BC}$ in $G$, we define the following vertex sets:
\begin{itemize}
  \item $P_{A}(G_{BC}) \defeq V_{mP}(G_{BC}, C_{G_{BC}} \cap A_{G})$ 
  \item $P_{NA}(G_{BC}) \defeq C_{G_{BC}} \setminus P_{A}(G_{BC})$
\end{itemize}

For $a \in ( C_{G_{BC}} \cap A_{G} )$,
\begin{itemize}
   \item $P'_{A}(G_{BC}, a) \defeq V_{mP}(G_{BC}, (C_{G_{BC}} \cap A_{G}) \setminus \{ a \})$ 
   \item $P'_{NA}(G_{BC}, a) \defeq C_{G_{BC}} \setminus P'_{A}(G_{BC}, a)$   
\end{itemize}

Let $R = C_{G_{BC}} \cap B_{G_{BC}}$ 
and define:
\begin{itemize}
    \item $R_{M} \defeq \{ r_t \in R \mid |P'_{NA}(G_{BC}, r_t)| = \max(\{ |P'_{NA}(G_{BC}, r)| \mid r \in R \}) \}$
\end{itemize}

For $r_{M} \in R_{M}$,
\begin{itemize}
    \item $
    P'_{NA}(G_{BC}) \defeq 
    \begin{cases}
      P'_{NA}(G_{BC}, r_{M}) & |R_{M}| > 0 \\
      P_{NA}(G_{BC}) & |R_{M}| = 0 \\
    \end{cases}
  $
\end{itemize}

Finally, define:

\begin{itemize}
    \item $P'_{A}(G_{BC}) \defeq C_{G_{BC}} \setminus P'_{NA}(G_{BC})$
\end{itemize}

\begin{figure}[tb]
    \begin{center}
    \includegraphics[width=0.6\linewidth]{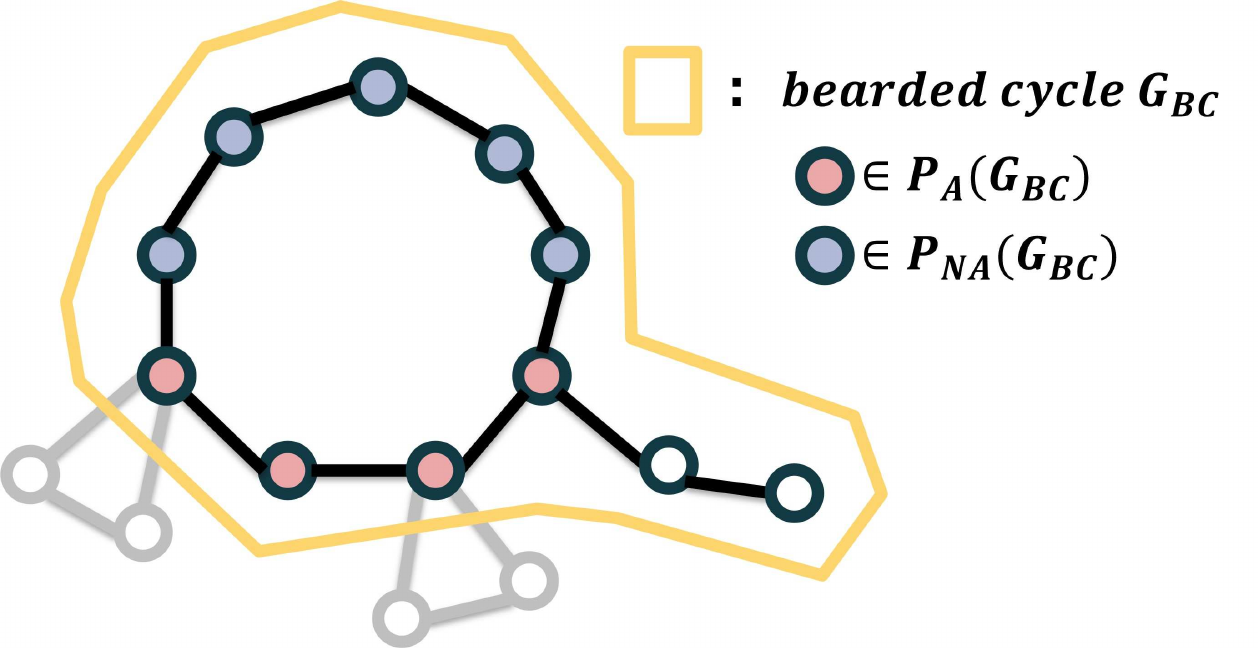}
    \caption{Examples of $P_{A}(G_{BC})$ and $P_{NA}(G_{BC})$ in a bearded cycle $G_{S}$}\label{fig:exp_papna} 
    \end{center}
\end{figure}

\begin{figure}[tb]
    \begin{center}
    \includegraphics[width=0.6\linewidth]{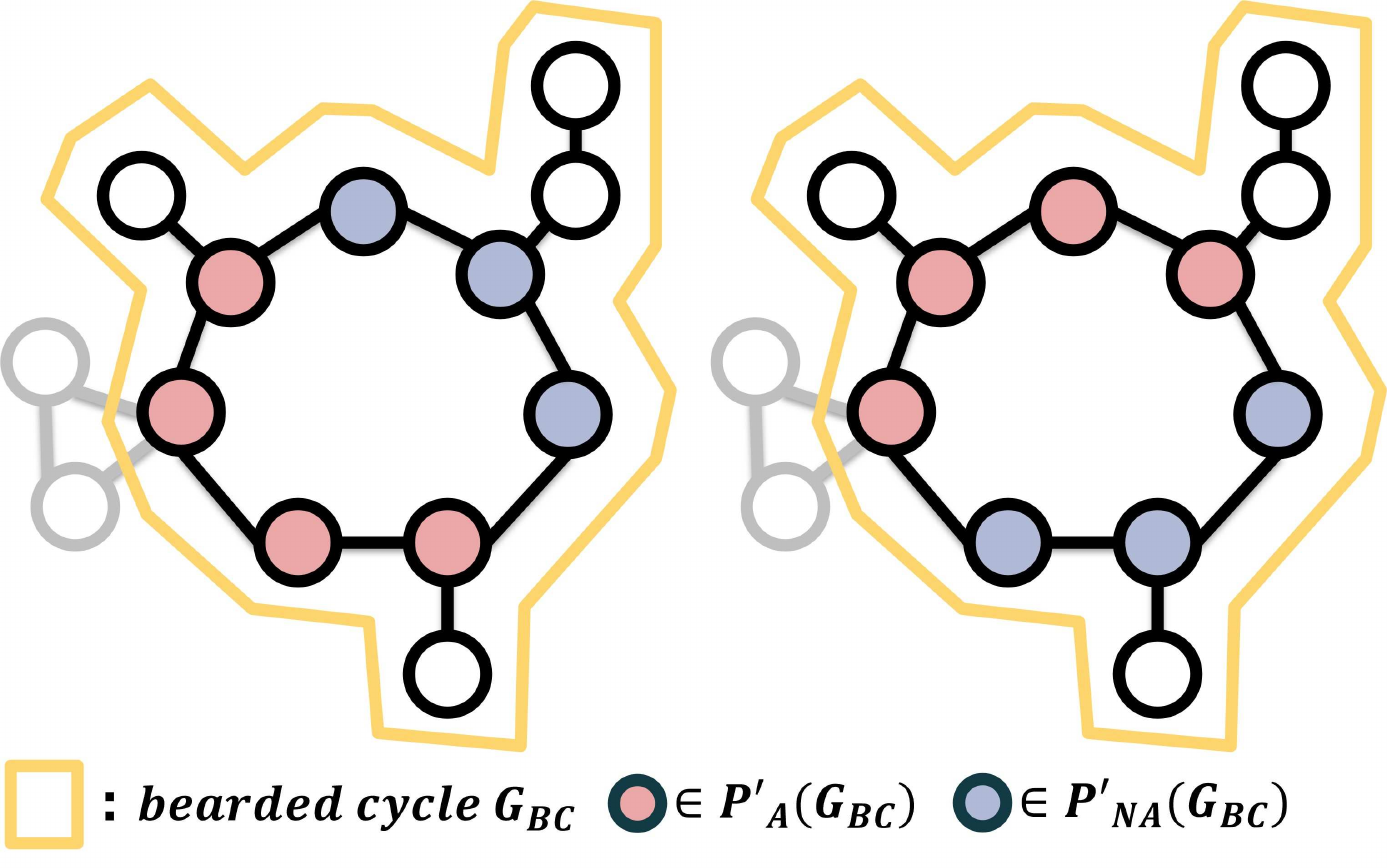}
    \caption{Examples of $P'_{A}(G_{BC})$ and $P'_{NA}(G_{BC})$ in a bearded cycle $G_{S}$}\label{fig:exp_pdapdna}
    \end{center}
\end{figure}

Examples of $P_{A}(G_{BC})$, $P_{NA}(G_{BC})$, $P'_{A}(G_{BC})$, and $P'_{NA}(G_{BC})$
are shown in Figures \rm{\ref{fig:exp_papna}} and \ref{fig:exp_pdapdna}.

Next, we define the concept of a $\subGPS$ for each \bc~and \nt~in a cactus graph $G$. 
The union of all $\subGPS$ sets over these components forms a $\GPS$ of the cactus graph. 
We compute the $\GPN$ by maximizing the size of the $\subGPS$ in each component.

\begin{definition} \rm{($\subGPS$)}
Let $G_S = (V_S, E_S)$ be a \bc~or \nt~in a cactus graph $G$. 
A vertex set $S_S \subseteq V_S$ is called a $\subGPS$ of $G_S$ if it satisfies all of the following:
\begin{itemize}
  \item For any $s_1, s_2 \in S_S$, no shortest path between $s_1$ and $s_2$ in $G_S$ contains any $s_3 \in S_S \setminus \{s_1, s_2\}$ (\ie $S_S$ is a $\GPS$ in $G_S$).
  \item For any $s_1 \in S_S$ and $a_1 \in (V_S \cap A_G) \setminus B_G$, no shortest path between $s_1$ and $a_1$ in $G_S$ contains any $s_2 \in S_S \setminus \{s_1\}$.
  \item For any $a_1, a_2 \in (V_S \cap A_G) \setminus B_G$, no shortest path between them in $G_S$ contains any $s_1 \in S_S$.
\end{itemize}
\end{definition}

It is worthwhile to note that the $\GPN$ of a cycle graph is \textbf{3} unless the number of vertices is 4.
As an exceptional case, the $\GPN$ of a cycle graph with 4 vertices is \textbf{2}.

Here we present the main theorem for computing the $\GPN$ of cactus graphs 
that are not simple cycles.

\begin{theorem}
\label{thm:main}
Let $G$ be a cactus graph that is not a simple cycle. \\
Let $BC_{all} = \{ BC_0, BC_1, \ldots, BC_k \}$ be the set of all \bc s in $G$. 
Then, the general position number $\GPN(G)$ is:
\begin{equation}
\notag
    \GPN(G) = 2n_1 + n_2 + n_3 + |L_G|
\end{equation}
where,
\begin{itemize}
  \item $n_1=|\{ BC\in BC_{all}\mid |P_{NA}(BC)|>|P_{A}(BC)|\}|$
  \item $n_2=|\{ BC\in BC_{all}\mid 0<|P_{NA}(BC)|\le |P_{A}(BC)|\}|$
  \item $n_3=|\{ BC\in BC_{all}\mid |P_{NA}(BC)|=0 \wedge |P'_{NA}(BC)|>|P'_{A}(BC)|\}|$
\end{itemize}
\end{theorem}


\section{The General Position Number of Cactus Graph}
In this section, we present the (exact) general position number ($\GPN$) of cactus graphs 
by proving Theorem \ref{thm:main} stated in the previous section. 
We begin with a lemma concerning $\subGPS$s in bearded cycles. 
This lemma allows us to estimate how many vertices on a cycle within a \bc~can be included in a $\subGPS$. 
%

Before the discussion about the size of each $\subGPS$,
we present two important lemmas about a special $\GPS$, named $\sGPS$,
which is a restricted $\GPS$ on a cycle graph when a subset of vertices $W \subseteq V$ is given.

\begin{definition} 
\rm{($\sGPS$)}
Given a cycle graph $G=(V,E)$ and a subset of vertices $W \subseteq V$,
a subset $S \subseteq V$ is called an $\sGPS$ if it satisfies all of the following:
\begin{itemize}
    \item For any pair of vertices $s_1, s_2 \in S$, no shortest path between them in $G$ contains any vertex $s_3 \in S$ such that $s_3 \ne s_1$ and $s_3 \ne s_2$.
    \item For any vertex $s_1 \in S$ and any vertex $w_1 \in W$, no shortest path between them in $G$ contains any vertex $s_2 \in S$ such that $s_2 \ne s_1$.
    \item For any pair of vertices $w_1, w_2 \in W$, no shortest path between them in $G$ contains any vertex $s_1 \in S$.
\end{itemize}
\end{definition}

Now we show two lemmas about the size of $\sGPS$ 
which are the keys for discussing the size of $\subGPS$.

\begin{lemma}
Let $G=(V,E)$ be a cycle graph, and let $W \subseteq V$ be a subset of the vertices. 
If $|W| > 0$ and $|P_{NW}| > |P_W|$, then the size of the maximum $\sGPS$ in $G$ is 2,
where $P_W = V_{mP}(G, W)$ and $P_{NW} = V \setminus P_W$. 
\label{lem:sgps2}
\end{lemma}

\begin{proof}
We first show that there exists an $\sGPS$ of size 2 in $G$.
Then, assuming that there exists an $\sGPS$ in $G$ of size greater than 2,
we prove by contradiction that this assumption leads to a contradiction.

Because $|P_{NW}| > |P_W|$, there exists an integer $d$ satisfying the following inequality:

\begin{equation}
\label{eq:range-d}
    0 \leq d < \frac{|P_{NW}|-|P_W|}{2}
\end{equation}

From this, the following inequality also holds:

\begin{equation}
\label{eq:PWPNW}
    |P_W|-1+d+1 < |P_{NW}|-1-d+1
\end{equation}

By the definitions of $P_W$ and $P_{NW}$,
the subgraph induced by $P_W$ and $P_{NW}$ are both connected.
Therefore, there must exist exactly two edges in $G$ that connect vectices in $P_W$ and $P_{NW}$.
Consider two distinct vertices $n_1, n_2 \in P_{NW}$ that are adjacent to vertices in $P_W$ 
(these exist by the definition of $P_{NW}$).
Let the vertices in $P_W$ adjacent to $n_1$ and $n_2$ be $w_1$ and $w_2$, respectively
(\ie $w_1, w_2 \in P_W$ and $(w_1, n_1), (w_2, n_2) \in E$).

We show that there exists a vertex set $\{s_1, s_2\} \subseteq P_{NW}$ satisfying the definition of an $\sGPS$.
Specifically, there exists a simple path from $w_1$ to $s_1$ in $G$ that does not include $w_2$,
of length $e+1$ ($e > 0$), and a simple path from $w_2$ to $s_2$ that does not include $w_1$,
also of length $e+1$.
Such paths must exist because $|P_{NW}| > |P_W| > 0$ (recall $|W|>0$).

If $e=d$, then from inequality (2), any shortest path from $s_1$ to $w_2$ does not include $s_2$,
and any shortest path from $s_2$ to $w_1$ does not include $s_1$.
Thus, $\{s_1, s_2\} \subseteq P_{NW}$ satisfies the definition of ans $\sGPS$.

To show maximality,
let $S = \{s_1, s_2\}$ and assume, for contradiction, that there exists another $\sGPS$ $S' \subseteq V$
such that $|S| < |S'|$.
By the definition of $\sGPS$, $S'$ must be a subset of $P_{NW}$.
Since $|W|>0$ and $|S'|>|S|=2$, there must exist some $s_1' \in S'$, 
and corresponding $s_2' \in S'$ and $w \in W$, 
such that a shortest path from $s_1'$ to $w$ includes $s_2'$.
This contradicts the definition of an $\sGPS$.
Therefore, such an $S'$ does not exist, and $S$ is a maximum $\sGPS$.
\end{proof}

\begin{figure}[tb] 
    \begin{center}
    \includegraphics[width=0.5\linewidth]{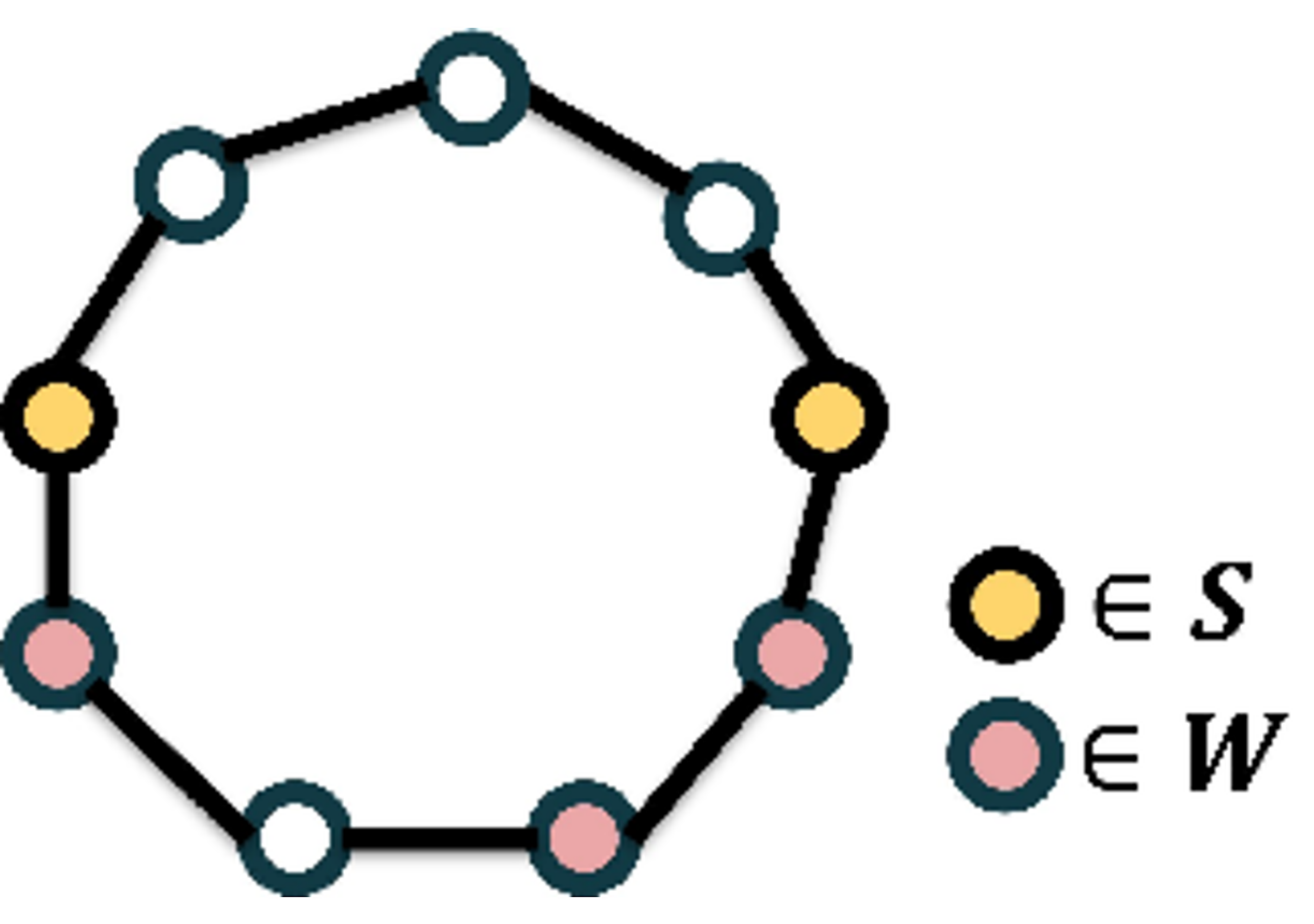}
    \caption{An example of $\sGPS$ $S$ with size 2 on a cycle graph $G$ when vertices set $W$ is given}
    \label{fig:exp_s-gps}
    \end{center}
\end{figure}

Figure \ref{fig:exp_s-gps} illustrates an example of $\sGPS$ of size 2.

The following lemma provides another case such that the size of the maximum $\sGPS$ is equal to 1.

\begin{lemma}
Let $G=(V,E)$ be a cycle graph, and let $W \subseteq V$ be a subset of the vertices. 
If $|W| > 0$ and $|P_{NW}| \le |P_W|$, then the size of the maximum $\sGPS$ in $G$ is equal to $1$, 
where $P_W = V_{mP}(G, W)$ and $P_{NW} = V \setminus P_W$. 
\label{lem:sgps1}
\end{lemma}

\begin{proof}
We first show that there exists an $\sGPS$ of size $1$ in $G$.
Then, assuming that there exists an $\sGPS$ in $G$ of size greater than $1$,
we prove by contradiction that this assumption leads to a contradiction.

To show existence, take any single vertex ${s}\subseteq P_{NW}$ (such an $s$ exists because $|P_{NW}|>0$).
By the definition of $P_W$, no shortest path between any $w_1,w_2\in P_W$ contains $s$.
Hence ${s}$ satisfies the definition of an $\sGPS$.

Let $S=\{s\}$. Assume, for the sake of contradiction, that there exists an $\sGPS$ $S'\subseteq V$ with $|S|<|S'|$.
We will obtain a contradiction by considering cases according to the size of $S'$.

\begin{enumerate}[\arabic*.]

\item
Case $|S'|=2$ : 
Write $S'=\{s'1,s'2\}$. By the definitions of $\sGPS$ and $P_W$, we have $S'\subseteq P_{NW}$.
From the definitions of $P_W$ and $P_{NW}$, the subgraphs induced by $P_W$ and $P_{NW}$ are each connected,
so there are two vertices in $P_{NW}$ adjacent to vertices in $P_W$.
Hence, because $0<|P_{NW}|$ and $|P_W|>0$, we have $2\le |P_{NW}|$ and therefore there exist two distinct vertices
$n_1,n_2\in P_{NW}$ that are adjacent to vertices in $P_W$.
Let $w_1,w_2\in P_W$ be vertices adjacent to $n_1$ and $n_2$, respectively (so $(w_1,n_1),(w_2,n_2)\in E$).
Let $d_1$ be the length of a shortest path from $s'_1$ to $n_1$, and $d_2$ the length of a shortest path from $s'_2$ to $n_2$.
Suppose that every shortest path from $s'1$ to $w_1$ does not contain $s'2$ or $w_2$.
Then, by the definition of $\sGPS$, no shortest path from $s'1$ to $w_2$ contains $s'2$, which yields
\begin{equation}
(d_1+1+|P_W|-1)-(|P_{NW}|-1-d_1-d_2+d_2+1)=2d_1+|P_W|-|P_{NW}|<0.
\tag{3}
\end{equation}
On the other hand, since $|P_{NW}|\le |P_W|$ and $d_1\ge 0$, we have
\begin{equation}
2d_1+|P_W|-|P_{NW}|\ge 0.
\tag{4}
\end{equation}
(3) contradicts (4). Therefore the assumption $|S'|=2$ leads to a contradiction.
\\
\item
Case $|S'|\ge 3$ : 
There must exist some $s'_1\in S'$ and some $s'_2\in S'$ and $w\in W$ such that a shortest path from $s'_1$ to $w$ contains $s'_2$.
This directly contradicts the definition of an $\sGPS$.

\end{enumerate}

In both cases we obtain a contradiction; hence no $\sGPS\ S'\subseteq V$ with $|S'|>|S|$ exists.
Therefore $S=\{s\}$ is a maximum $\sGPS$ of $G$.
\end{proof}


Now we present important lemmas that provides the size of the maximum $\subGPS$ 
in each type of \bc s and \nt s. First, we show the lemma about $\subGPS$.

\begin{lemma}
For a \bc~or \nt~$G_S=(V_S,E_S)$ in a cactus graph $G$ that is not a cycle graph, for every $i$ such that $1\le i \le $ the size of maximum $\subGPS$ in $G_S$, there exists a $\subGPS$ of size $i$ that contains no articulation point($\in A_G$).
\label{lem:subGPSnoArticulation}
\end{lemma}
\begin{proof}
Let $S$ be any $\subGPS$ in $G_S$ and suppose, by contradiction, that there exists a vertex $a\in S\cap A_G$. Let $(V_{1},E_{1}),(V_{2},E_{2}),...,(V_{k},E_{k})\,(k\ge 2)$ be the new connected components of $G_S-\{a\}$, created by removing $a$. 
By the definition of $\subGPS$, there is only one index $i\in \{1,...,k\}$ such that $S\cap V_i \neq \emptyset$, otherwise there would be two vertices $u,v$ belonging to two different connected components in $G_S-\{a\}$ that $\{u,v\}$ is a $\subGPS$ in $G_S$. This is impossible since any shortest path between $u$ and $v$ passes through $a$. 
Then we have these property below.
There exists some vertex  $u\in V_{j}\,(i\neq j\wedge V_J\cap ((V_S\cap A_G)\setminus B_G)=\emptyset)$ such that any shortest path between $u$ and any vertex in $S\cap V_{i}$ containing any other vertex of $S\cap V_{i}$.
Additionally, $u$ is not contained in any shortest path between any two vertices in $(V_S\cap A_G)\setminus B_G$(If $u$ is contained, then there exist some shortest paths between $(V_S\cap A_G)\setminus B_G$ that containing $a$ twice, contradicting the definition of a shortest path). 
Moreover, any shortest path between $u$ and any vertex in $(V_S\cap A_G)\setminus B_G$ is not containing any vertex in $S\cap V_{i}$.
Moreover, by the definition of $\subGPS$, $a$ is not contained in any shortest path between any two vertices in $S\setminus \{ a\}$.
From these above, then $S'=(S\setminus \{a \})\cup \{ u\}$ is a $\subGPS$ in $G_S$.
\end{proof}

Next, we show the lemmas that provides the size of the maximum $\subGPS$ of \bc~in each type.

\begin{lemma}
Let $G$ be a cactus graph that is not a cycle graph, and let $G_S = (V_S, E_S)$ be a \bc~
satisfying $|P_{NA}(G_S)| > |P_A(G_S)|$ in $G$. 
The size of the maximum $\subGPS$ in $G_S$ is equal to $2 + |L_{G_S}|$.
\label{lem:cond1}
\end{lemma}
\begin{proof}
We first show that there exists a $\subGPS$ of size $2 + |L_{G_S}|$ in $G_S$. 
We then show that another $\subGPS$ whose size is greater than $2 + |L_{G_S}|$ does not exist.

To prove the existence of a $\subGPS$ of size $2 + |L_{G_S}|$, note that by the definition of $\subGPS$ and from $(C_{G_S} \cap A_G) \setminus B_{G_S} \subseteq C_{G_S} \cap A_G$ and Lemma \ref{lem:sgps2}, there exists a $\subGPS$ $\{s_1, s_2\} \subseteq P_{NA}(G_S)$. Moreover, since $L_{G_S}$ is the set of leaves, it satisfies the conditions of $\subGPS$ (leaves are not contained in any shortest path between any two other vertices).

From $B_{G_S} \cap C_{G_S} \subseteq P_A(G_S)$, it follows that no shortest path between any two vertices in $L_{G_S}$ contains any vertex in $P_{NA}(G_S)$. Thus, $\{s_1, s_2\} \cup L_{G_S}$ satisfies the definition of a $\subGPS$. Let $S = \{s_1, s_2\} \cup L_{G_S}$.

Suppose, for contradiction, that there exists a $\subGPS$ $S' \subseteq V_S$ with $|S| < |S'|$ and $S'\cap A_G=\emptyset$(from lemma \ref{lem:subGPSnoArticulation}, there exists such a $\subGPS$). We divide the discussion based on the cases defined by $B_{G_S}$ and $(C_{G_S} \cap A_G) \setminus B_{G_S}$ (i.e., the existence of beards in $G_S$, and the existence of articulation points in $G$ that are in $G_S$ but not in a beard). Before discussion, since $|S'|>2+|L_{G_S}|$ and $(V_S \setminus C_{G_S}) \cap S' \subseteq L_{G_S}$(because $S' \cap  A_G=\emptyset$), we have $|C_{G_S} \cap S'| > 2$.

\begin{enumerate}[\arabic*.]

\item 
Case $B_{G_S} = \emptyset$ and $(C_{G_S} \cap A_G) \setminus B_{G_S} \ne \emptyset$ : 
Since $L_{G_S} = \emptyset$, $|S| = 2$. From $(C_{G_S} \cap A_G) \setminus B_{G_S} \ne \emptyset$ and Lemma \ref{lem:sgps2}, the maximum $\subGPS$ size in $G_S$ is 2. This contradicts the assumption that $|S'| > |S|$.
\\
\item 
Case $B_{G_S} \ne \emptyset$ and $(C_{G_S} \cap A_G) \setminus B_{G_S} = \emptyset$ : 
If $|C_{G_S} \cap S'|> 3$, then $C_{G_S} \cap S'$ violates the $\subGPS$ condition. If $|C_{G_S} \cap S'| = 3$, then in any beard $B'$, there exists some $c_1 \in C_{G_S} \cap S'$ such that a shortest path from every vertex in $B'$ to $c_1$ includes some $c_2 \in C_{G_S} \cap S'$ ($c_2 \ne c_1$), implying that $(V_S \setminus C_{G_S}) \cap S' = \emptyset$, and thus $|S'| = 3$. Since $B_{G_S} \ne \emptyset$, we have $|S| \ge 3$, contradicting $|S'| > |S|$.
\\
\item 
Case $B_{G_S} \ne \emptyset$ and $(C_{G_S} \cap A_G) \setminus B_{G_S} \ne \emptyset$ : 
From $|C_{G_S} \cap S'| > 2$ and $(C_{G_S} \cap A_G) \setminus B_{G_S} \ne \emptyset$ and Lemma \ref{lem:sgps2}, we have $|C_{G_S} \cap S'| \le 2$, leading to a contradiction.

\end{enumerate}

Hence, assuming the existence of a $\subGPS$ $S' \subseteq V_S$ with $|S| < |S'|$ leads to a contradiction. Therefore, $S$ is a maximum $\subGPS$ in $G_S$.
\end{proof}

\begin{figure}[tbhp]
        \begin{minipage}[htpb]{0.49\columnwidth}
            \begin{center}
            \includegraphics[width=0.8\linewidth]{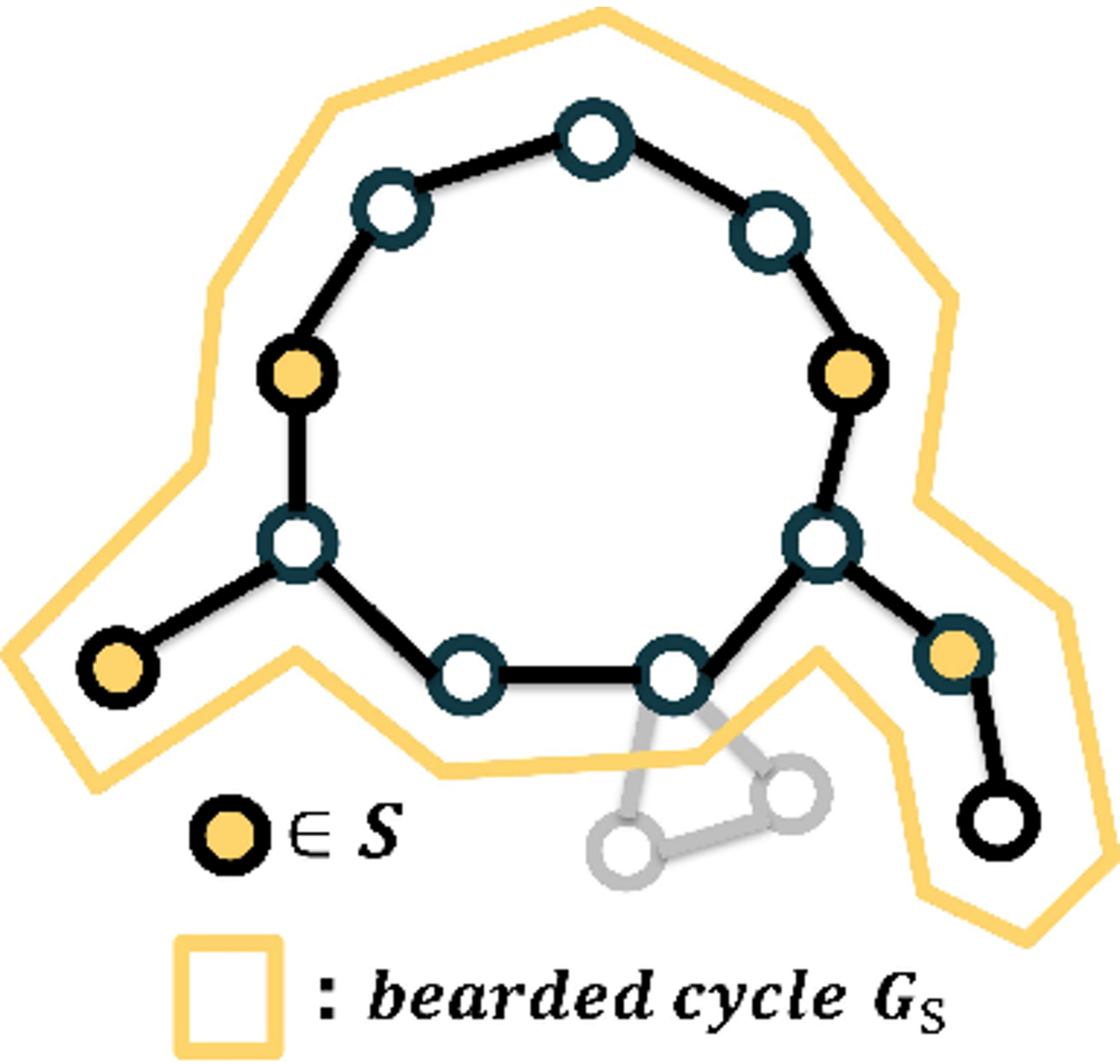}
            \caption{A \bc~satisfying \textbf{Cond.1}}\label{fig:beardcycle_subGPS_condition1}
            \end{center}
        \end{minipage}   
        \hfill
        \begin{minipage}[htbp]{0.49\columnwidth}
            \begin{center}
            \includegraphics[width=0.8\linewidth]{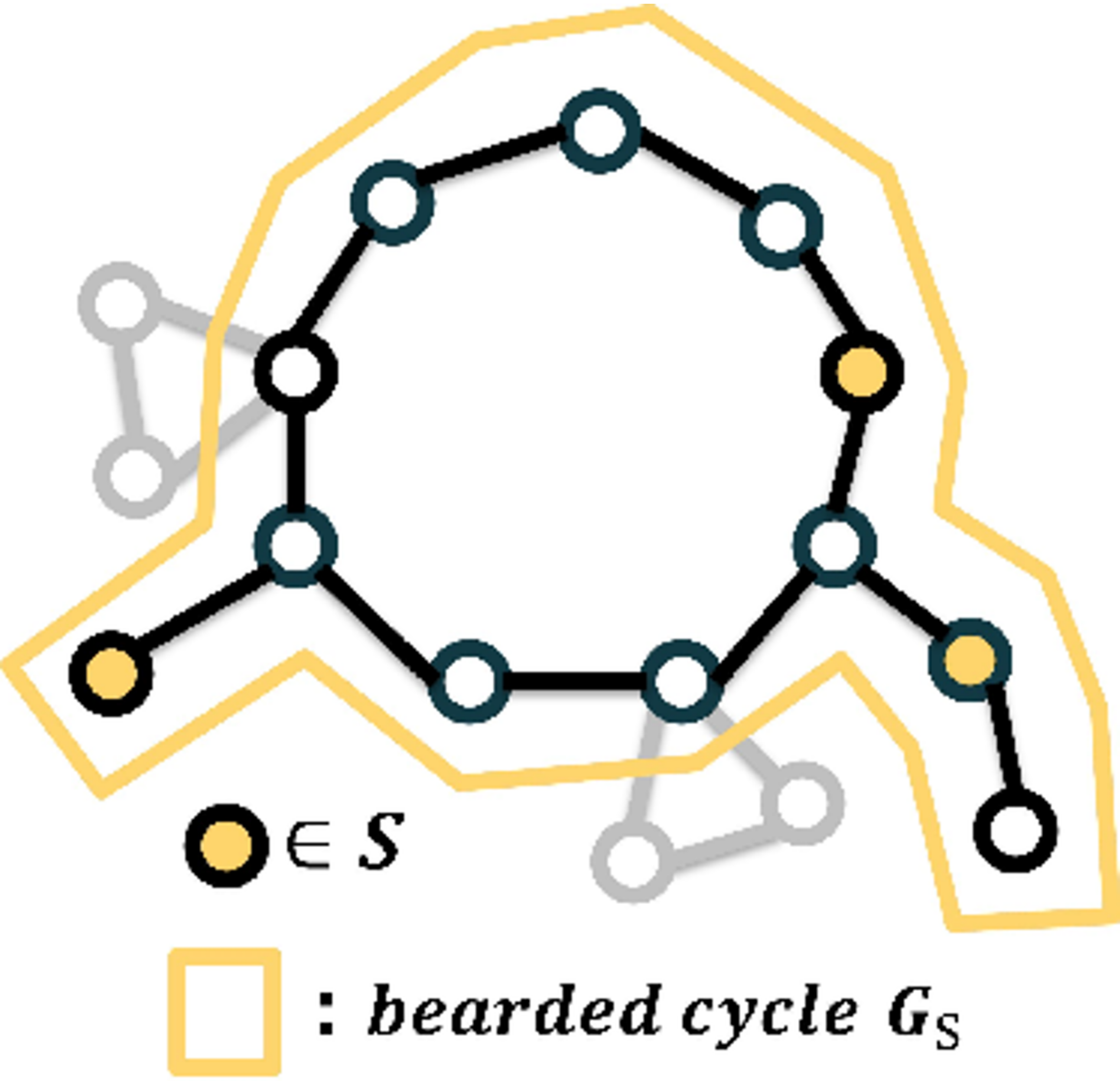}
            \caption{A \bc~satisfying \textbf{Cond.2}}\label{fig:beardcycle_subGPS_condition2}
            \end{center}
        \end{minipage}

        \vspace{20pt}
        \begin{minipage}[htpb]{0.49\columnwidth}
            \begin{center}
            \includegraphics[width=0.81\linewidth]{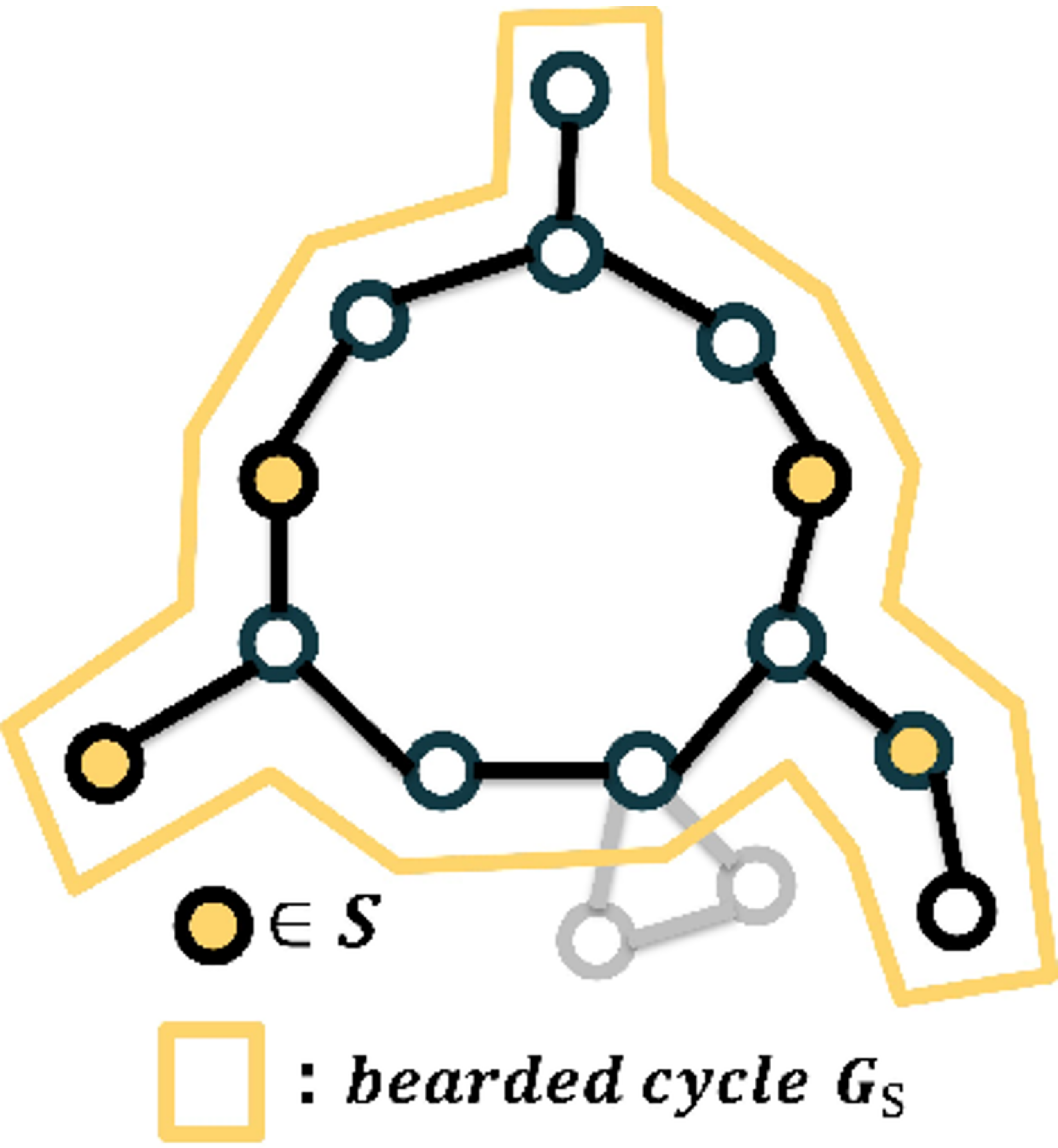}
            \caption{A \bc~satisfying \textbf{Cond.3}}\label{fig:beardcycle_subGPS_condition3}
            \end{center}
        \end{minipage}   
        \hfill
        \begin{minipage}[htbp]{0.49\columnwidth}
            \begin{center}
            \includegraphics[width=0.8\linewidth]{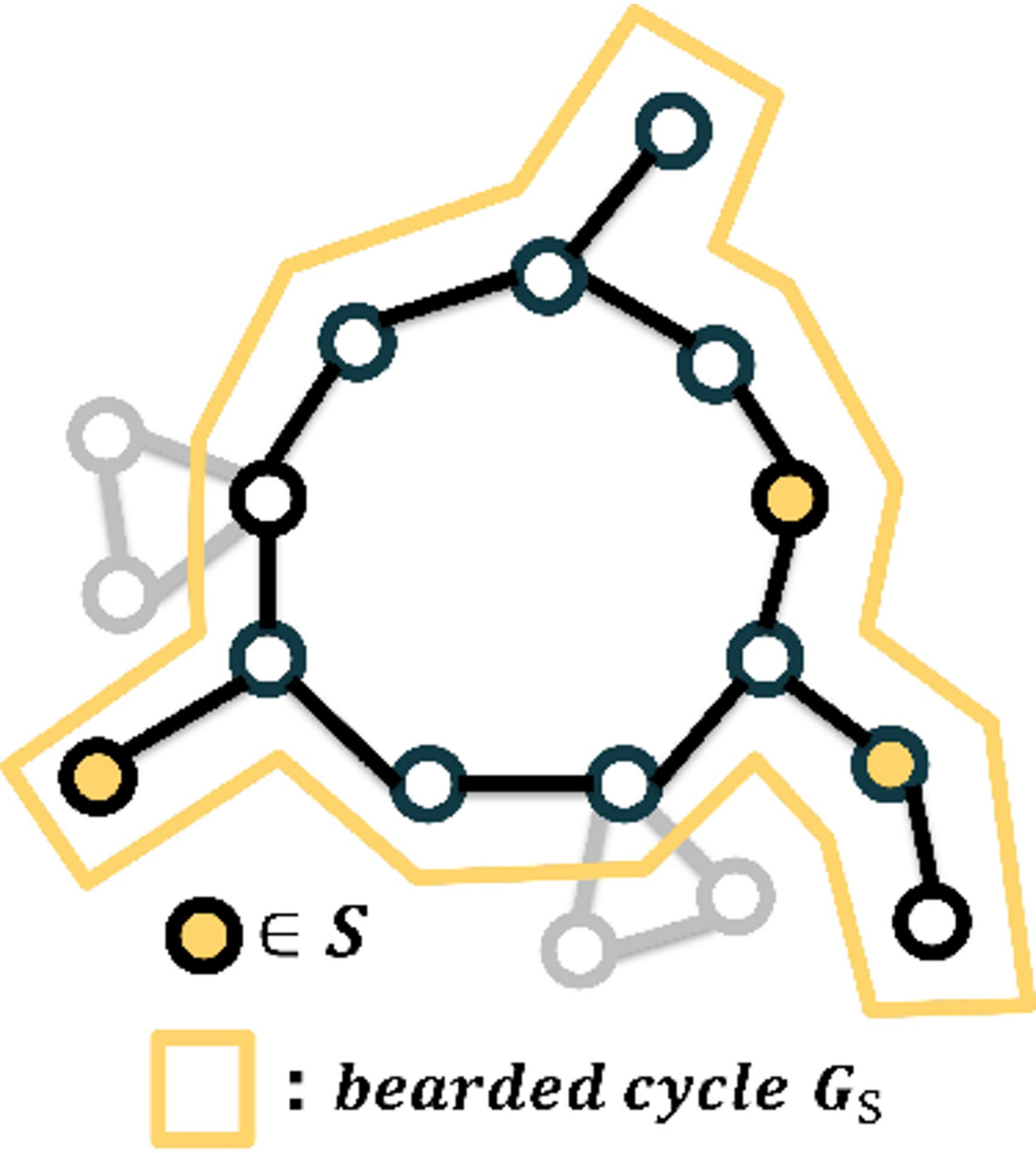}
            \caption{A \bc~satisfying \textbf{Cond.4}}\label{fig:beardcycle_subGPS_condition4}
            \end{center}
        \end{minipage}

        \vspace{20pt}
        \begin{minipage}[htpb]{0.49\columnwidth}
            \begin{center}
            \includegraphics[width=0.81\linewidth]{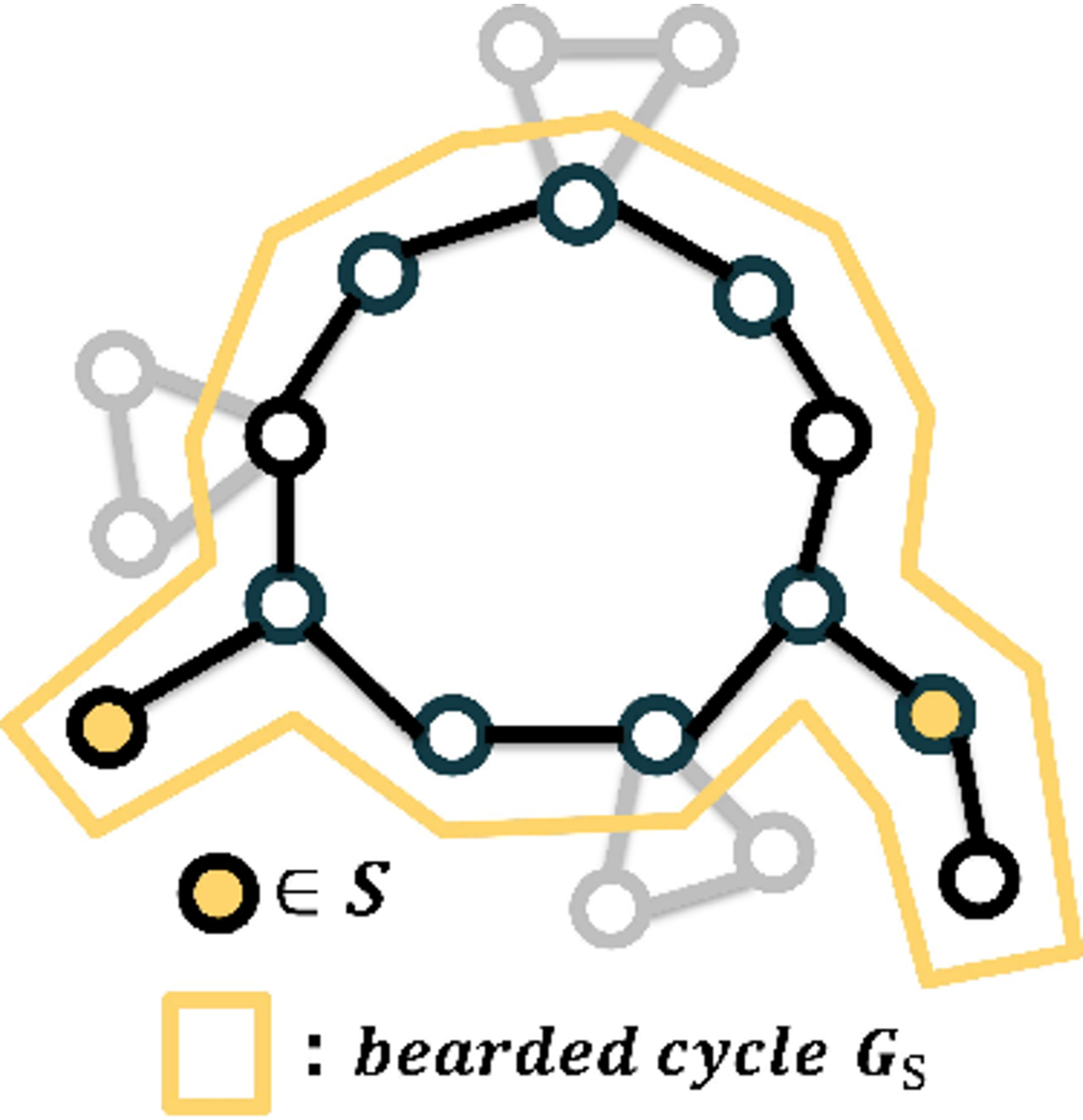}
            \caption{A \bc~satisfying \textbf{Cond.5}}\label{fig:beardcycle_subGPS_condition5}
            \end{center}
        \end{minipage}   
        \hfill
        \begin{minipage}[htbp]{0.49\columnwidth}
            \begin{center}
            \includegraphics[width=0.8\linewidth]{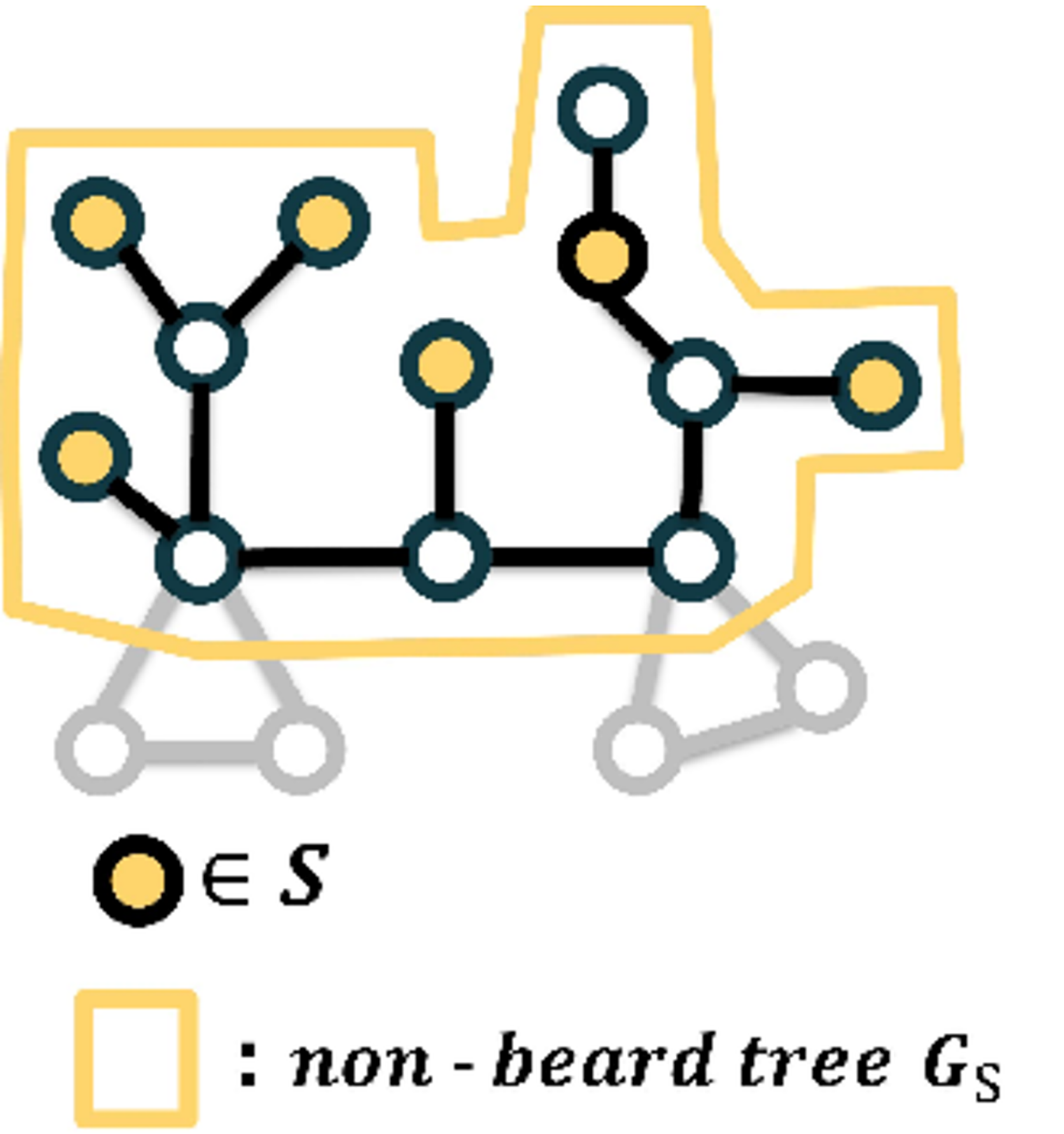}
            \caption{Example of a \nt}\label{fig:unbeardtree_subGPS}
            \end{center}
        \end{minipage}
\end{figure}

Figure \ref{fig:beardcycle_subGPS_condition1} illustrates an example of $G_S$ and maximum $\subGPS$ $S$.

\begin{lemma}
Let $G$ be a cactus graph that is not a cycle graph, and let $G_S = (V_S, E_S)$ be a \bc~in $G$. If $0<|P_{NA}(G_S)| \le |P_A(G_S)|$, then the size of the maximum $\subGPS$ in $G_S$ is equal to $1 + |L_{G_S}|$.
\label{lem:cond2}
\end{lemma}

\begin{proof}
We first show that there exists a $\subGPS$ of size $1 + |L_{G_S}|$ in $G_S$.
We then show that another $\subGPS$ whose size is greater than $1 + |L_{G_S}|$ does not exist.

To prove the existence of a $\subGPS$ of size $1 + |L_{G_S}|$, note that by the definition of $\subGPS$ and from $(C_{G_S}\cap A_G)\setminus B_{G_S}\subseteq C_{G_S}\cap A_G$, and Lemma \ref{lem:sgps1}, any one vertex $\{s\}\subseteq P_{NA}(G_S)$ satisfies the definition of a $\subGPS$.
Moreover, since $L_{G_S}$ is the set of leaves, every vertex of $L_{G_S}$ also satisfies the conditions of $\subGPS$ (leaves are not contained in any shortest path between any two other vertices).
By the definitions of $P_A(G_S)$ and $P_{NA}(G_S)$, for any $s\in P_{NA}(G_S)$,  no shortest path between any two vertices of $L_{G_S}$ contains $s$.
Therefore ${s}\cup L_{G_S}$ satisfies the definition of a $\subGPS$.
Let $S={s}\cup L_{G_S}$.

Now suppose, for contradiction, that there exists a $\subGPS$ $S'\subseteq V_S$ with $|S|<|S'|$ and $S'\cap A_G=\emptyset$(from lemma \ref{lem:subGPSnoArticulation}, there exists such a $\subGPS$).
We consider cases according to $B_{G_S}$ and $(C_{G_S}\cap A_G)\setminus B_{G_S}$. Before discussion, since $|S'|>1+|L_{G_S}|$ and $(V_S \setminus C_{G_S}) \cap S' \subseteq L_{G_S}$(because $S' \cap  A_G=\emptyset$), we have $|C_{G_S} \cap S'| > 1$.

\begin{enumerate}[\arabic*.]

\item 
Case $B_{G_S}=\emptyset$ and $(C_{G_S}\cap A_G)\setminus B_{G_S}\neq\emptyset$ : 
In this case $|L_{G_S}|=0$, hence $|S|=1$. By $(C_{G_S}\cap A_G)\setminus B_{G_S}\neq\emptyset$ and Lemma \ref{lem:sgps1}, the maximum size of a $\subGPS$ in $G_S$ is $1$. This contradicts the assumption that $|S'|>|S|$.

\item 
Case $B_{G_S}\neq\emptyset$ and $(C_{G_S}\cap A_G)\setminus B_{G_S}=\emptyset$ : 
Note that $2\le |C_{G_S}\cap S'|\le 3$ (if $|C_{G_S}\cap S'|\ge 4$ then $C_{G_S}\cap S'$ violates the $\subGPS$ condition).
We split further by the value of $|C_{G_S}\cap S'|$.

\begin{enumerate}[\alph*.]
     
\item
Case $|C_{G_S}\cap S'|=3$ : 
By an argument analogous to that in (ii) of the proof of Lemma \ref{lem:cond1}, we obtain $|S'|=3$.
To compare $|S'|$ with $|S|$ we estimate $|L_{G_S}|$. Since $0<|P_{NA}(G_S)|\le |P_A(G_S)|$, we have
\begin{gather}
2<|C_{G_S}|=|P_{NA}(G_S)|+|P_A(G_S)|\le 2|P_A(G_S)|,\
1<|P_A(G_S)|,\
1<|L_{G_S}|.
\end{gather}
Hence $|S|=1+|L_{G_S}|\ge 3$, so $|S'|=3$ contradicts $|S'|>|S|$.

\item
Case $|C_{G_S}\cap S'|=2$ : 
Consider two distinct vertices $u_1,u_2\in P_{NA}(G_S)$ that are adjacent to vertices in $P_A(G_S)$.
Let $v_1,v_2\in P_A(G_S)$ be vertices adjacent to $u_1,u_2$, respectively ($(u_1,v_1),(u_2,v_2)\in E_S$).
Compare simple path lengths among $u_1,u_2,v_1,v_2$.
A path from $u_2$ to $v_1$ not passing through $v_2$ has length $|P_{NA}(G_S)|-1+1$, while a path from $u_2$ to $v_1$ that passes through $v_2$ has length $1+|P_A(G_S)|-1$.
Since $0<|P_{NA}(G_S)|\le |P_A(G_S)|$, we obtain
\begin{gather}
(|P_{NA}(G_S)|-1+1)-(1+|P_A(G_S)|-1)=|P_{NA}(G_S)|-|P_A(G_S)|\le 0.
\end{gather}
Therefore every vertex of $P_{NA}(G_S)$ lies on every shortest path from $u_2$ to $v_1$.
Consequently, for any two vertices of $P_{NA}(G_S)$, one of them lies on every shortest path from the other to $v_1$.
We now consider how $S'$ meets the beards.

\begin{enumerate}[\roman*.]
    
\item
Case Each beard of $G_S$ contains exactly one vertex of $S'$ : 
Then no vertex of $C_{G_S}\cap S'$ belongs to $V_{mP}(G_S,(V_S\setminus C_{G_S})\cap S')$, and by $(C_{G_S}\cap A_G)\setminus B_{G_S}=\emptyset$ and the definition of $P_A(G_S)$, we have $C_{G_S}\cap S'\subseteq P_{NA}(G_S)$.
There exists a $s'1\in C_{G_S}\cap S'$ and $s'2\in (V_S\setminus C_{G_S})\cap S'$ such that every shortest path from $s'1$ to $s'2$ contains another vertex of $C{G_S}\cap S'$ different from $s'_1$(By the observation above (every vertex of $P_{NA}(G_S)$ lies on every shortest path from $u_2$ to $v_1$)). 
This contradicts the $\subGPS$ condition.

\item
Case There exists a beard of $G_S$ that contains no vertex of $S'$ : 
Then $|S'|\le 2+|L_{G_S}|-1=1+|L_{G_S}|$, contradicting the assumption $|S'|>|S|=1+|L_{G_S}|$.
\end{enumerate}
\end{enumerate}

\item
Case $B_{G_S}\neq\emptyset$ and $(C_{G_S}\cap A_G)\setminus B_{G_S}\neq\emptyset$ : 
As before $2\le |C_{G_S}\cap S'|\le 3$(if $4\le |C_{G_S}\cap S'|$, it violates the $\subGPS$ condition). We proceed by splitting into subcases according to the value of $|C_{G_S}\cap S'|$.

\begin{enumerate}[\alph*.]
    
\item
Case $|C_{G_S}\cap S'|=3$ : 
Since $(C_{G_S}\cap A_G)\setminus B_{G_S}\neq\emptyset$, and the definitions of $P_A(G_S)$ and $\subGPS$, and Lemmas \ref{lem:sgps2} and \ref{lem:sgps1}, we obtain $|C_{G_S}\cap S'|\le 2$. This Contradicts the $|C_{G_S}\cap S'|=3$

\item
Case $|C_{G_S}\cap S'|=2$ : 
If every beard contains exactly one vertex of $S'$, then by the same reasoning as in 2.b.i. we find $C_{G_S}\cap S'\subseteq P_{NA}(G_S)$ and, using the same path-length argument as before,
there exists $s'1\in C_{G_S}\cap S'$ and $s'2\in (V_S\setminus C{G_S})\cap S'$ such that every shortest path from $s'1$ to $s'2$ contains another vertex of $C_{G_S}\cap S'$, contradicting the $\subGPS$ condition.
If some beard contains no vertex of $S'$, then $|S'|\le 2+|L_{G_S}|-1=1+|L_{G_S}|$, again contradicting $|S'|>|S|$.
\end{enumerate}
\end{enumerate}

In all cases the assumption that a $\subGPS$ $S'\subseteq V_S$ with $|S'|>|S|$ and $S'\cap A_G=\emptyset$ exists leads to a contradiction.
Therefore no such $S'$ exists, and $S$ is a maximum $\subGPS$ in $G_S$.
\end{proof}

\begin{lemma}
Let $G$ be a cactus graph that is not a cycle graph, and let $G_S = (V_S, E_S)$ be a \bc~in $G$. If $|P_{NA}(G_{S})| =0 \wedge |P'_{NA}(G_{S})|>|P'_{A}(G_{S})|$, then the size of the maximum $\subGPS$ in $G_S$ is equal to $1 + |L_{G_S}|$.
\label{lem:cond3}
\end{lemma}

\begin{proof}
First we show that there exists a $\subGPS$ of size $1+|L_{G_S}|$ in $G_S$. We then show that another $\subGPS$ whose size is greater than $1 + |L_{G_S}|$ does not exist.

We first construct a $\subGPS$ $S$ of size $1+|L_{G_S}|$.
By the definition of $P'_{A}(G_S)$, we have $V_{mP}(G_{S},(V_{S}\cap A_{G})\setminus B_{G_{S}})\subseteq P'_{A}(G_S)$. 
Using this and definition of $\subGPS$, $|P'_{NA}(G_S)|>|P'_{A}(G_S)|$, Lemma \ref{lem:sgps2}, we obtain a set of two vertices $\{s_1,s_2\}\subseteq P'_{NA}(G_S)$ that satisfy the $\subGPS$ condition.
Next consider vertices in $(V_S\setminus C_{G_S})\cap S$. Let $B'_{G_S}$ denote the union of all beards that share at least one vertex with $P'_{A}(G_S)$. Since $B'_{G_S}\cap L_{G_S}$ is a set of leaves, it satisfies the $\subGPS$ condition.
Next, we claim that $\{s_1,s_2\}\cup (B'_{G_S}\cap L_{G_S})$ is a $\subGPS$.
$\{s_1,s_2\}\subseteq P'_{NA}(G_S)$ and $|P'_{NA}(G_S)|>|P'_{A}(G_S)|$ imply that for each vertex of $B'_{G_S}\cap L_{G_S}$, each shortest path from $s_1$ to that leaf does not contain $s_2$, and symmetrically each shortest path from $s_2$ to that leaf does not contain $s_1$. Moreover, by the definition of $P'_{A}(G_S)$ no vertex of $P'_{NA}(G_S)$ lies on any shortest path between two vertices of $B'_{G_S}\cap L_{G_S}$. Hence $\{s_1,s_2\}\cup (B'_{G_S}\cap L_{G_S})$ satisfies the $\subGPS$ condition. Let $S=\{ s_{1}, s_{2}\}\cup (B'_{G_{S}}\cap L_{G_{S}})$.
We will show that $S$ is maximum by deriving a contradiction from the assumption that there exists a $\subGPS$ $S'\subseteq V_S$ with $|S'|>|S|$ and $S'\cap A_G=\emptyset$(from lemma \ref{lem:subGPSnoArticulation}, there exists such a $\subGPS$).

Note that $|P_{NA}(G_S)|=0$ and $|P'_{NA}(G_S)|>|P'_{A}(G_S)|$ imply $B_{G_S}\neq\emptyset$. We therefore split into cases depending on whether $(C_{G_S}\cap A_G)\setminus B_{G_S}$ is empty or not. 
Additionally, because $S' \cap  A_G=\emptyset$, we have $(V_S \setminus C_{G_S}) \cap S' \subseteq L_{G_S}$($|(V_S \setminus C_{G_S}) \cap S'| \le |L_{G_S}|$).
From $|S'|\ge 2+|L_{G_S}|$ and $|(V_S \setminus C_{G_S}) \cap S'| \le |L_{G_S}|$ it follows that $|C_{G_{S}}\cap S'|\ge 2$(if $|(V_S \setminus C_{G_S}) \cap S'| = |L_{G_S}|$, from $|P_{NA}(G_S)|=0$, we have $|S'|=|L_{G_S}|<|S|$, a contradiction).

\begin{enumerate}[\arabic*.]

\item 
Case $(C_{G_S}\cap A_G)\setminus B_{G_S}\neq\emptyset$ : 
Because $(C_{G_S}\cap A_G)\setminus B_{G_S}\neq\emptyset$ and by Lemmas \ref{lem:sgps2} and \ref{lem:sgps1}, we have $|C_{G_S}\cap S'|\le 2$, hence $|C_{G_S}\cap S'|=2$.
Since $|S'|\ge 2+|L_{G_S}|$ and $|C_{G_S}\cap S'|=2$ and $|(V_S \setminus C_{G_S}) \cap S'| \le |L_{G_S}|$, $|(V_S \setminus C_{G_S}) \cap S'| = |L_{G_S}|$. 
From this and $|P_{NA}(G_S)|=0$, $|C_{G_S}\cap S'|=0$, a contradiction.

\item 
Case $(C_{G_S}\cap A_G)\setminus B_{G_S}=\emptyset$ : 
Since $(C_{G_S}\cap A_G)\setminus B_{G_S}=\emptyset$, we have $2\le |C_{G_S}\cap S'|\le 3$ (if $|C_{G_S}\cap S'|\ge 4$ then $C_{G_S}\cap S'$ would contradicts the $\subGPS$ condition). We proceed by splitting into subcases according to the value of $|C_{G_S}\cap S'|$.

\begin{enumerate}[\alph*.]

\item
Case $|C_{G_S}\cap S'|=2$ : 
Since $|S'|\ge 2+|L_{G_S}|$ and $|C_{G_S}\cap S'|=2$ and $|(V_S \setminus C_{G_S}) \cap S'| \le |L_{G_S}|$, $|(V_S \setminus C_{G_S}) \cap S'| = |L_{G_S}|$. 
From this and $|P_{NA}(G_S)|=0$, $|C_{G_S}\cap S'|=0$, a contradiction.

\item
Case $|C_{G_S}\cap S'|=3$ : 
For any leaf $b\in L_{G_S}$, there exists $c\in C_{G_S}\cap S'$ such that some shortest path from $b$ to $c$ contains another vertex of $C_{G_S}\cap S'$ different from $c$. Therefore $(B_{G_S}\setminus C_{G_S})\cap S'=\emptyset$, and thus $|S'|=|C_{G_S}\cap S'|=3$.
To compare with $|S|$, note that $(C_{G_S}\cap A_G)\setminus B_{G_S}=\emptyset$ and $|P_{NA}(G_S)|=0$ imply the existence of two vertices $r_1,r_2\in C_{G_S}\cap B_{G_S}$ such that every vertex of $C_{G_S}$ lies on some shortest path between $r_1$ and $r_2$. From this it follows that $|L_{G_S}|\ge 2$, hence $|S|=1+|L_{G_S}|\ge 3$. Therefore $|S'|=3$  contradicts $|S'|>|S|$.
\end{enumerate}
\end{enumerate}

In all possible cases the assumption that a $\subGPS$ $S'\subseteq V_S$ with $|S'|>|S|$ and $S'\cap A_G=\emptyset$ exists leads to a contradiction. Therefore no such $S'$ exists, and $S$ is a maximum $\subGPS$ in $G_S$.

\end{proof}

\begin{lemma}
Let $G$ be a cactus graph that is not a cycle graph, and let $G_S = (V_S, E_S)$ be a \bc~in $G$. If $|P_{NA}(G_{S})| =0 \wedge 0<|P'_{NA}(G_{S})|\le |P'_{A}(G_{S})|$, then the size of the maximum $\subGPS$ in $G_S$ is equal to $|L_{G_S}|$.
\label{lem:cond4}
\end{lemma}

\begin{proof}
We first show that there exists a $\subGPS$ of size $|L_{G_S}|$ in $G_S$. 
We then show that another $\subGPS$ whose size is greater than $|L_{G_S}|$ does not exist.

Since every leaf vertex does not lie on any shortest path between two other vertices, the set $L_{G_S}$ satisfies the definition of a $\subGPS$. Let $S=L_{G_S}$.

Suppose there exists a $\subGPS$ $S'\subseteq V_S$ with $|S'|>|S|=|L_{G_S}|$ and $S'\cap A_G=\emptyset$(from lemma \ref{lem:subGPSnoArticulation}, there exists such a $\subGPS$). 
From this we have $(V_S \setminus C_{G_S}) \cap S' \subseteq L_{G_S}$($|(V_S \setminus C_{G_S}) \cap S'| \le |L_{G_S}|$).
From $|S'|\ge 1+|L_{G_S}|$ and $|(V_S \setminus C_{G_S}) \cap S'| \le |L_{G_S}|$ it follows that $|C_{G_{S}}\cap S'|\ge 1$(if $|(V_S \setminus C_{G_S}) \cap S'| = |L_{G_S}|$, from $|P_{NA}(G_S)|=0$, we have $|S|=|S'|=|L_{G_S}|$, a contradiction).

By the definition of $\subGPS$, $1\le |C_{G_S}\cap S'|\le 3$ (if $|C_{G_S}\cap S'|\ge 4$ then $C_{G_S}\cap S'$ would contradict the definition of $\subGPS$). We split according to the value of $|C_{G_S}\cap S'|$.

\begin{enumerate}[\arabic*.]

\item
Case $|C_{G_S}\cap S'|=1$ : 
Then, from $|(V_S \setminus C_{G_S}) \cap S'| \le |L_{G_S}|$ and $|S'|\ge 1+|L_{G_S}|$, we obtain $|(B_{G_{S}}\setminus C_{G_{S}})\cap S'|= |L_{G_S}|$. 
Since every beard contains exactly one vertex of $S'$, $P_A(G_S)=C_{G_S}$ and the definition of $\subGPS$ imply $|C_{G_S}\cap S'|=0$, contradicting $|C_{G_S}\cap S'|=1$.

\item
Case $|C_{G_S}\cap S'|=2$ : 
From $|(V_S \setminus C_{G_S}) \cap S'| \le |L_{G_S}|$ and $|S'|\ge 1+|L_{G_S}|$, we have $|L_{G_S}|-1\le|(B_{G_{S}}\setminus C_{G_{S}})\cap S'| \le |L_{G_S}|$. 

If $|(B_{G_{S}}\setminus C_{G_{S}})\cap S'| = |L_{G_S}|$, then every beard contains exactly one vertex of $S'$, $P_A(G_S)=C_{G_S}$ and the definition of $\subGPS$ force $|C_{G_S}\cap S'|=0$, contradicting $|C_{G_S}\cap S'|=2$.

Consider the remaining possibility: $|(B_{G_S}\setminus C_{G_S})\cap S'|=|L_{G_S}|-1$. In this case there is exactly one beard that contains no vertex of $S'$. Let $B'_E$ denote the set of all vertices belonging to that beard. By the definition of $P'_A(G_S)$ we have $|P'_{A}(G_S)|\le |V_{mP}(G_{S},(C_{G_{S}}\cap A_{G})\setminus B'_{E})|$. Sinse $|P'_{NA}(G_S)|\le |P'_A(G_S)|$ and by Lemma \ref{lem:sgps1}, it follows that $|C_{G_S}\cap S'|\le 1$, contradicting $|C_{G_S}\cap S'|=2$.

\item
Case $|C_{G_S}\cap S'|=3$ : 
First assume $(C_{G_S}\cap A_G)\setminus B_{G_S}\neq\emptyset$. Then there exist vertices $c\in C_{G_S}\cap S'$ and $a\in (C_{G_S}\cap A_G)\setminus B_{G_S}$ such that some shortest path from $c$ to $a$ contains another vertex of $C_{G_S}\cap S'$, contradicting the definition of $\subGPS$.
Next assume $(C_{G_S}\cap A_G)\setminus B_{G_S}=\emptyset$. Then for every vertex $b\in B_{G_S}\setminus C_{G_S}$ there exists $c\in C_{G_S}\cap S'$ such that some shortest path from $b$ to $c$ contains another vertex of $C_{G_S}\cap S'$. Hence $(B_{G_S}\setminus C_{G_S})\cap S'=\emptyset$ and therefore $|S'|=|C_{G_S}\cap S'|=3$.
We now estimate $|S|$. From $1\le |P'_{NA}(G_{S})|\le |P'_A(G_{S})| \wedge 3\le |C_{G_{S}}|=|P'_{NA}(G_{S})|+|P'_{NA}(G_{S})|$, we deduce $1\le |P'_{NA}(G_S)| \wedge 2\le |P'_A(G_S)|$. Using $(C_{G_S}\cap A_G)\setminus B_{G_S}=\emptyset$ and the definition of $P'_{NA}(G_S)$, it follows that $|L_{G_S}|\ge 3$. Hence $|S|=|L_{G_S}|\ge 3$, so $|S'|=3$ does not exceed $|S|$, contradicting $|S'|>|S|$.

\end{enumerate}

In all possible cases, the assumption that there exists a $\subGPS$ $S'\subseteq V_S$ with $|S'|>|S|$ and $S'\cap A_G=\emptyset$ leads to a contradiction. Therefore no such $S'$ exists, and $S=L_{G_S}$ is a maximum $\subGPS$ in $G_S$.
\end{proof}

\begin{lemma}
Let $G$ be a cactus graph that is not a cycle graph, and let $G_S = (V_S, E_S)$ be a \bc~in $G$. If $|P'_{NA}(G_{S})| =0$, then the size of the maximum $\subGPS$ in $G_S$ is equal to $|L_{G_S}|$.
\label{lem:cond5}
\end{lemma}

\begin{proof}
We first show that there exists a $\subGPS$ of size $|L_{G_S}|$ in $G_S$. 
We then show that another $\subGPS$ whose size is greater than $|L_{G_S}|$ does not exist.

Let $S=L_{G_S}$. Since every leaf vertex does not lie on any shortest path between two other vertices, $S$ satisfies the definition of a $\subGPS$.

Suppose, for contradiction, that there exists a $\subGPS$ $S' \subseteq V_S$ with $|S| < |S'|$ and $S'\cap A_G=\emptyset$(from lemma \ref{lem:subGPSnoArticulation}, there exists such a $\subGPS$). 
From this, we have $(V_S \setminus C_{G_S}) \cap S' \subseteq L_{G_S}$($|(V_S \setminus C_{G_S}) \cap S'| \le |L_{G_S}|$). From this and $|S'|\ge 1+|L_{G_S}|$ it follows that $|C_{G_{S}}\cap S'|\ge 1$(if $|(V_S \setminus C_{G_S}) \cap S'| = |L_{G_S}|$, from $|P'_{NA}(G_S)|=0$, we have $|S'|=|L_{G_S}|<|S|$, a contradiction).

We divide the discussion based on the cases defined by $B_{G_S}$ and $(C_{G_S} \cap A_G) \setminus B_{G_S}$.

\begin{enumerate}[\arabic*.]

\item 
Case $B_{G_S}\neq\emptyset\wedge (C_{G_S}\cap A_G)\setminus B_{G_S}\neq\emptyset$ : 
Then, from $1\le |C_{G_S}\cap S'|$ together with Lemmas \ref{lem:sgps2} and\ref{lem:sgps1} we have $1\le |C_{G_S}\cap S'|\le 2$. We consider the possible values.

\begin{enumerate}[\alph*.]

\item
Case $|C_{G_S}\cap S'|=1$ : 
Since $|S'|\ge 1+|L_{G_S}|$, we have $|(B_{G_{S}}\setminus C_{G_{S}})\cap S'|= |L_{G_{S}}|$. 
Since every beard contains exactly one vertex of $S'$, $P_A(G_S)=C_{G_S}$ and the definition of a $\subGPS$ forces $|C_{G_S}\cap S'|=0$, contradicting $|C_{G_S}\cap S'|=1$.

\item
Case $|C_{G_S}\cap S'|=2$ : 
From $|S'|\ge 1+|L_{G_S}|$, we get $|L_{G_{S}}|-1\le|(B_{G_{S}}\setminus C_{G_{S}})\cap S'|\le |L_{G_{S}}|$.
If every beard contains exactly one vertex of $S'$, $P_A(G_S)=C_{G_S}$ and the definition of a $\subGPS$ forces $|C_{G_S}\cap S'|=0$, contradicting $|C_{G_S}\cap S'|=2$.
There remains the possibility that exactly one beard contains no vertex of $S'$. Let $B'_E$ denote the set of vertices of that beard. By the definition of $P'_A(G_S)$, we have $|P'_{A}(G_{S})|\le |V_{mP}(G_{S},(C_{G_{S}}\cap A_{G})\setminus B'_{E})|$. Since $P'_A(G_S)=C_{G_S}$($|P'_{NA}(G_S)|=0$), it follows that $|C_{G_S}\cap S'|=0$, contradicting $|C_{G_S}\cap S'|=2$.

\end{enumerate}

\item 
Case $B_{G_S}=\emptyset \wedge (C_{G_S}\cap A_G)\setminus B_{G_S}\neq\emptyset$ : 
Here $|P'_{NA}(G_S)|=0$ and $B{G_S}=\emptyset$, by the definition of a $\subGPS$, we obtain $|C_{G_S}\cap S'|=0$, contradicting $|C_{G_S}\cap S'|\ge 1$.

\item 
Case $B_{G_S}\neq\emptyset \wedge (C_{G_S}\cap A_G)\setminus B_{G_S}=\emptyset$ : 
In this case, $1\le |C_{G_S}\cap S'|\le 3$. We divide the discussion based on $|C_{G_S}\cap S'|$.

\begin{enumerate}[\alph*.]

\item
Case $|C_{G_S}\cap S'|=1$ : 
Since $|S'|\ge 1+|L_{G_S}|$, we have $|(B_{G_S}\setminus C_{G_S})\cap S'|= |L_{G_S}|$. 
Since every beard contains exactly one vertex of $S'$, $P_A(G_S)=C_{G_S}$ forces $|C_{G_S}\cap S'|=0$, contradiction.

\item
Case $|C_{G_S}\cap S'|=2$ : 
From $|S'|\ge 1+|L_{G_S}|$, we get $|L_{G_{S}}|-1\le|(B_{G_{S}}\setminus C_{G_{S}})\cap S'|\le |L_{G_{S}}|$.
If every beard contains exactly one vertex of $S'$, $P_A(G_S)=C_{G_S}$ forces $|C_{G_S}\cap S'|=0$, contradiction.
There remains the possibility that exactly one beard contains no vertex of $S'$. Let $B'_E$ denote the set of vertices of that beard. By the definition of $P'_A(G_S)$, we have $|P'_{A}(G_{S})|\le |V_{mP}(G_{S},(C_{G_{S}}\cap A_{G})\setminus B'_{E})|$. Since $P'_A(G_S)=C_{G_S}$($|P'_{NA}(G_S)|=0$), it follows that $|C_{G_S}\cap S'|=0$, contradicting $|C_{G_S}\cap S'|=2$.

\item
Case $|C_{G_S}\cap S'|=3$ : 
For every $b\in B_{G_S}\setminus C_{G_S}$ there exists $c\in C_{G_S}\cap S'$ such that some shortest path from $b$ to $c$ contains another vertex of $C_{G_S}\cap S'$. Hence $(B_{G_S}\setminus C_{G_S})\cap S'=\emptyset$ and therefore $|S'|=3$. 
Moreover, from $B_{G_S}\neq\emptyset$, $(C_{G_S}\cap A_G)\setminus B_{G_S}=\emptyset$, and the definition of $P'_A(G_S)$, there exists vertices $r_1,r_2\in B_{G_S}\cap C_{G_S}$ such that every vertex of $C_{G_S}$ lies on some shortest path between $r_1$ and $r_2$. This implies $|L_{G_S}|\ge 3$, hence $|S|=|L_{G_S}|\ge 3$, so $|S'|=3$ does not exceed $|S|$, contradicting $|S'|>|S|$.

\end{enumerate}
\end{enumerate}

In all possible cases the assumption that a $\subGPS$ $S'\subseteq V_S$ with $|S'|>|S|$ and $S'\cap A_G=\emptyset$ exists leads to a contradiction. Therefore no such $S'$ exists, and $S=L_{G_S}$ is a maximum $\subGPS$ in $G_S$.
\end{proof}

\begin{lemma}
Let $G$ be a cactus graph that is not a cycle graph, and let $G_S = (V_S, E_S)$ be a \nt~in $G$, then the size of the maximum $\subGPS$ in $G_S$ is equal to $|L_{G_S}|$.
\label{lem:nonbt}
\end{lemma}

\begin{proof}
We first show that there exists a $\subGPS$ of size $|L_{G_S}|$ in $G_S$. 
We then show that another $\subGPS$ whose size is greater than $|L_{G_S}|$ does not exist.

Each leaf vertex is not contained in any shortest path between any two other vertices, therefore $L_{G_S}$ satisfies the definition of a $\subGPS$. Let $S=L_{G_S}$. 

Suppose, for contradiction, that there exists a $\subGPS$ $S' \subseteq V_S$ with $|S| < |S'|$ and $S'\cap A_G=\emptyset$(from lemma \ref{lem:subGPSnoArticulation}, there exists such a $\subGPS$). 

From $|S'|>|S|=|L_{G_S}|$, we have $S\cap A_G\neq\emptyset$. This contradicts $S'\cap A_G=\emptyset$. 
Therefore no such $S'$ exists, and $S$ is a maximum $\subGPS$ in $G_S$.
\end{proof}

From Lemma \ref{lem:cond1} to Lemma \ref{lem:nonbt},
$\subGPS$ of every component, \bc s and \nt s, in $G$ can be determined.
Figures \ref{fig:beardcycle_subGPS_condition1} to \ref{fig:unbeardtree_subGPS} illustrate 
each example of maximum $\subGPS$s for each type of \bc~and \nt.

The following lemma (Lemma \ref{lem:unionsubGPS}) shows that the union of all $\subGPS$s becomes a legitimate $\GPS$ of a given cactus graph $G$,
and the next lemma (Lemma \ref{lem:sizeofunion}) provides the size of $\GPS$ consisting of all $\subGPS$s.

\begin{lemma}
A union of $\subGPS$s of all \bc s and \nt s in a non-cycle cactus graph $G$
becomes a $\GPS$ of $G$.
\label{lem:unionsubGPS}
\end{lemma}

\begin{proof}
Let two \bc s or \nt s of $G$ be $G_{S_1}=(V_{S_1},E_{S_1})$, $G_{S_2}=(V_{S_2},E_{S_2})$ with $G_{S_1}\neq G_{S_2}$. Let $S_1$ and $S_2$ be maximum $\subGPS$s of $G_{S_1}$ and $G_{S_2}$, respectively.
We prove that $S_1 \cup S_2$ is a $\GPS$ of $G$ by checking the following argument.

\begin{enumerate}[\arabic*.]

\item

We show that $S_1$ is a $GPS$ in $G$.
By the definition of $\subGPS$, for any $s_1,s_2\in S_1$, no shortest path between $s_1$ and $s_2$ in $G_{S_1}$ contains any vertex $s_3 \in S_1 (s_3 \neq s_1,s_2)$. Since shortest paths in $G_{S_1}$ are also shortest paths in $G$. Hence, by the definition of a $\GPS$, $S_1$ is a $\GPS$ of $G$.

\item
We show that for any $s_1\in S_1$ and any $a\in (V_{S_1}\cap A_G)\setminus B_G$, no shortest path from $s_1$ to $a$ contains any $s_2\in S_1$ with $s_2\neq s_1$.
From the definition of a $\subGPS$ in $G_{S_1}$, for any $s_1\in S_1$ and any $a\in (V_{S_1}\cap A_G)\setminus B_G$, no shortest path from $s_1$ to $a$ contains any $s_2\in S_1$ with $s_2\neq s_1$ in $G$.

\item
We show that for any $a_1,a_2\in (V_{S_1}\cap A_G)\setminus B_G$, no shortest path between $a_1$ and $a_2$ contains any vertex $s\in S_1$.
From the definition of a $\subGPS$ in $G_{S_1}$, for any $a_1,a_2\in (V_{S_1}\cap A_G)\setminus B_G$, no shortest path between $a_1$ and $a_2$ contains any vertex $s\in S_1$ in $G$.

\end{enumerate}

From (i), $S_1$ and $S_2$ are $\GPS$s in $G$. If $V_{S_1}\cap V_{S_2}\neq \emptyset$, then, from (i) and (ii), $S_1 \cup S_2$ is a $\GPS$ of $G$. If $V_{S_1}\cap V_{S_2}= \emptyset$, then, from (i) and (ii) and (iii), $S_1 \cup S_2$ is a $\GPS$ of $G$.

Therefore, the union of the $\subGPS$s of the \bc s and the \nt s of $G$ is a $\GPS$ of $G$.
\end{proof}

\begin{lemma}
Let $G$ be a cactus graph that is not a simple cycle. \\
Let $BC_{all} = \{ BC_0, BC_1, \ldots, BC_k \}$ be the set of all \bc s in $G$. 
A size of union of $\subGPS$ of all \bc s and \nt s in $G$ is:
\begin{equation}
\notag
    2n_1 + n_2 + n_3 + |L_G|
\end{equation}
where,
\begin{itemize}
  \item $n_1=|\{ BC\in BC_{all}\mid |P_{NA}(BC)|>|P_{A}(BC)|\}|$
  \item $n_2=|\{ BC\in BC_{all}\mid 0<|P_{NA}(BC)|\le |P_{A}(BC)|\}|$
  \item $n_3=|\{ BC\in BC_{all}\mid |P_{NA}(BC)|=0 \wedge |P'_{NA}(BC)|>|P'_{A}(BC)|\}|$
\end{itemize}
\label{lem:sizeofunion}
\end{lemma}

\begin{proof}
We prove the claim by summing the sizes of maximum $\subGPS$s of all \bc s and all \nt s of $G$. 
By lemma \ref{lem:cond1}...\ref{lem:cond5}, define $p_1(BC)$, $p_2(BC)$, $p_3(BC)$, $p_4(BC)$, $p_5(BC)$ on a \bc $BC\in BC_{all}$ as follows :  
\begin{itemize}
  \item $p_{1}(BC) : |P_{NA}(BC)|>|P_{A}(BC)|$
  \item $p_{2}(BC) : 0<|P_{NA}(BC)|\le |P_{A}(BC)|$
  \item $p_{3}(BC) : |P_{NA}(BC)| =0 \wedge |P'_{NA}(BC)|>|P'_{A}(BC)|$
  \item $p_{4}(BC) : |P_{NA}(BC)| =0 \wedge 0<|P'_{NA}(BC)|\le |P'_{A}(BC)|$
  \item $p_{5}(BC) : |P_{NA}(BC)| =0 \wedge 0<|P'_{NA}(BC)|=0$
\end{itemize}
For $i=1,...,5$, define the integer $n_i=\{ BC\in BC_{all} \mid p_{i}(BC) \}$.
Let $UBT_{all}$ be the set of all \nt s in $G$. Let $S$ be the union of the maximum $\subGPS$s of all \bc s and all \nt s in $G$.
Then the size $|S|$ is given as follows.
\begin{equation}
  \begin{split}
  |S|=&\sum_{BC\in \{ BC\in BC_{all}\, \mid \,  p_{1}(BC) \}} (2+|L_{BC}|)+\sum_{BC\in \{ BC\in BC_{all}\, \mid \,  p_{2}(BC) \}} (1+|L_{BC}|) \\
     &+\sum_{BC\in \{ BC\in BC_{all}\, \mid \,  p_{3}(BC) \}} (1+|L_{BC}|)+\sum_{BC\in \{ BC\in BC_{all}\, \mid \,  p_{4}(BC) \}} (|L_{BC}|) \\
     &+\sum_{BC\in \{ BC\in BC_{all}\, \mid \,  p_{5}(BC) \}} (|L_{BC}|)+\sum_{UBT\in UBT_{all}} (|L_{UBT}|) \\
     =&\sum_{BC\in \{ BC\in BC_{all}\, \mid \,  p_{1}(BC) \}} (2)+\sum_{BC\in \{ BC\in BC_{all}\, \mid \,  p_{2}(BC) \}} (1) \\
     &+\sum_{BC\in \{ BC\in BC_{all}\, \mid \,  p_{3}(BC) \}} (1)+|L_{G}| \\
     =&2n_{1}+n_{2}+n_{3}+|L_{G}|
  \end{split}
\end{equation}
\end{proof}

The two following lemmas are necessary to prove the main theorem (Theorem \ref{thm:main}).

\begin{lemma}
Let $G$ be a cactus graph that is not a cycle graph, and let $G_S = (V_S, E_S)$ be a \bc~or \nt~in $G$, and $S$ be a union set of maximum $\subGPS$ in all \bc s and \nt s in $G$, then each connected component of $G-V_{S}\setminus (A_{G}\setminus B_{G_{S}})-E_{S}$ includes at least 2 vertices in $S$.
\label{lem:subGPSAroundGS}
\end{lemma}
\begin{proof}
Let $G-V_{S}\setminus (A_{G}\setminus B_{G_{S}})-E_{S}=G_{NS}=(V_{NS},E_{NS})$.
We prove the claim by considering each connected compornent of $G_{NS}$. 
Let $G_{NSC}=(V_{NSC},E_{NSC})$ be an arbitary connected compornent of $G_{NS}$. We divide the discussion based on the value of $|L_G \cap V_{NSC}|$.

\begin{enumerate}[\arabic*.]

\item
Case $|L_G \cap V_{NSC}|\ge 2$ : 
By Lemmas \ref{lem:cond1},...\ref{lem:nonbt}, we have $|V_{NSC}\cap S|\ge 2$.

\item
Case $|L_G \cap V_{NSC}|=1$ : 
We divide the discussion based on which subgraph contains the vertex of $L_G \cap V_{NSC}$

\begin{enumerate}[\alph*.]  

\item
Case A \nt contains the vertex of $L_G \cap V_{NSC}$ : 
First, $G_{NSC}$ must contain a \bc . If $G_{NSC}$ does not contain a \bc , then, by the definitions of \bc and \nt , the single leaf vertex in $L_{G}\cap V_{NSC}$ is not included in a \nt , contradicting a assumption. Hence $G_{NSC}$ contains some bearded cycle $G_{BC}$.
If no \bc in $G_{NSC}$ satisfies $|(C_{G_{BC}}\cap A_{G})\setminus B_{G_{BC}}|=1$, then contradicting the definitions of \bc~ and \nt , so there exists a \bc in $G_{NSC}$ satisfies $|(C_{G_{BC}}\cap A_{G})\setminus B_{G_{BC}}|=1$. 
Because this $G_{BC}$ contains no beards and 
$|(C_{G_{BC}}\cap A_{G})\setminus B_{G_{BC}}|=1$, we have $|C_{G_{BC}}\cap A_{G}|=1$, then $|P_{NA}(G_{BC})|>|P_{A}(G_{BC})|$, and by Lemma \ref{lem:cond1}, it follows that $|V_{BC}\cap S|=2$. Therefore $|V_{NSC}\cap S|\ge 2$.

\item
Case A \bc contains the vertex of $L_G \cap V_{NSC}$ : 
Let $G_{BC}$ be a bearded cycle that contains the vertex in $L_{G}\cap V_{NSC}$. 
If $|(C_{G_{BC}}\cap A_{G})\setminus B_{G_{BC}}|=1$, then $|P'_{NA}(G_{BC})|>|P'_{A}(G_{BC})|$, and by Lemma \ref{lem:cond3}, it follows that $|V_{BC}\cap S|=2$. Therefore $|V_{NSC}\cap S|\ge 2$.

If $|(C_{G_{BC}}\cap A_{G})\setminus B_{G_{BC}}|\ge 2$, then, by an argument in 2.a. , there exists another \bc~ $G'_{BC}=(V'_{BC},E'_{BC})$ in $G_{NSC}$ with $|(C_{G'_{BC}}\cap A_{G})\setminus B_{G'_{BC}}|=1$. Then  $|P_{NA}(G'_{BC})|>|P_{A}(G'_{BC})|$, it follows that $|V'_{BC}\cap S|=2$, so $|V_{NSC}\cap S|\ge 2$.

\end{enumerate}

\item
Case $|L_G\cap V_{NSC}|=0$ : 
By an argument in 2.a. , there exists a \bc~ $G'_{BC}=(V'_{BC},E'_{BC})$ in $G_{NSC}$ with $|(C_{G'_{BC}}\cap A_{G})\setminus B_{G'_{BC}}|=1$. Then  $|P_{NA}(G'_{BC})|>|P_{A}(G'_{BC})|$, it follows that $|V'_{BC}\cap S|=2$, so $|V_{NSC}\cap S|\ge 2$.

\end{enumerate}

In all cases, we have $|V_{NSC}\cap S|\ge 2$. Hence every connected compornent of $G-V_{S}\setminus (A_{G}\setminus B_{G_{S}})-E_{S}$ includes at least 2 vertices in $S$.
\end{proof}

\begin{lemma}
For a graph $G$, for every $i$ such that $1\le i \le \GPN(G)$, there exists a $\GPS$ of size $i$ that contains no articulation point.
\end{lemma}
\begin{proof}
Let $S$ be any $\GPS$ and suppose, by contradiction, that there exists a vertex $a\in S\cap A_G$. Let $(V_{1},E_{1}),(V_{2},E_{2}),...,(V_{k},E_{k})\,(k\ge 2)$ be the new connected components of $G-\{a\}$, created by removing $a$. By the definition of $\GPS$, there is only one index $i\in \{1,...,k\}$ such that $S\cap V_i \neq \emptyset$, otherwise there would be two vertices $u,v$ belonging to two different connected components in $G-\{a\}$ that $\{u,v\}$ is a $\GPS$ in $G$. This is impossible since any shortest path between $u$ and $v$ passes through $a$. Then $S'=(S\setminus \{a \})\cup \{ a'\}$, where $a'\in V_{j}\,(j\neq i)$, is such that $S'$ is a $\GPS$ in $G$.
\end{proof}

Next, we summarize the conditions for each type of \bc~$G_{S}=(V_{S},E_{S})$ 
and the size of the maximum $\subGPS$ for each type as the following:

\begin{itemize}
    \item \textbf{Condition 1}: $|P_{NA}(G_{S})|>|P_{A}(G_{S})|$\\ 
    \textbf{Maximum size of $\subGPS$}: $2+|V_{S}\cap L|$ (Lemma \ref{lem:cond1})
    \item \textbf{Condition 2}: $0<|P_{NA}(G_{S})|\le |P_{A}(G_{S})|$\\
    \textbf{Maximum size of $\subGPS$}: $1+|V_{S}\cap L|$ (Lemma \ref{lem:cond2})
    \item \textbf{Condition 3}: $|P_{NA}(G_{S})|=0\wedge |P'_{NA}(G_{S})|>|P'_{A}(G_{S})|$\\
    \textbf{Maximum size of $\subGPS$}: $1+|V_{S}\cap L|$ (Lemma \ref{lem:cond3})
    \item \textbf{Condition 4}: $|P_{NA}(G_{S})|=0\wedge 0<|P'_{NA}(G_{S})|\le |P'_{A}(G_{S})|$\\
    \textbf{Maximum size of $\subGPS$}: $|V_{S}\cap L|$ (Lemma \ref{lem:cond4})
    \item \textbf{Condition 5}: $|P_{NA}(G_{S})|=0\wedge |P'_{NA}(G_{S})|=0$\\
    \textbf{Maximum size of $\subGPS$}: $|V_{S}\cap L|$ (Lemma \ref{lem:cond5})
\end{itemize}

\begin{figure}[tb]
    \begin{center}
    \includegraphics[width=0.8\linewidth]{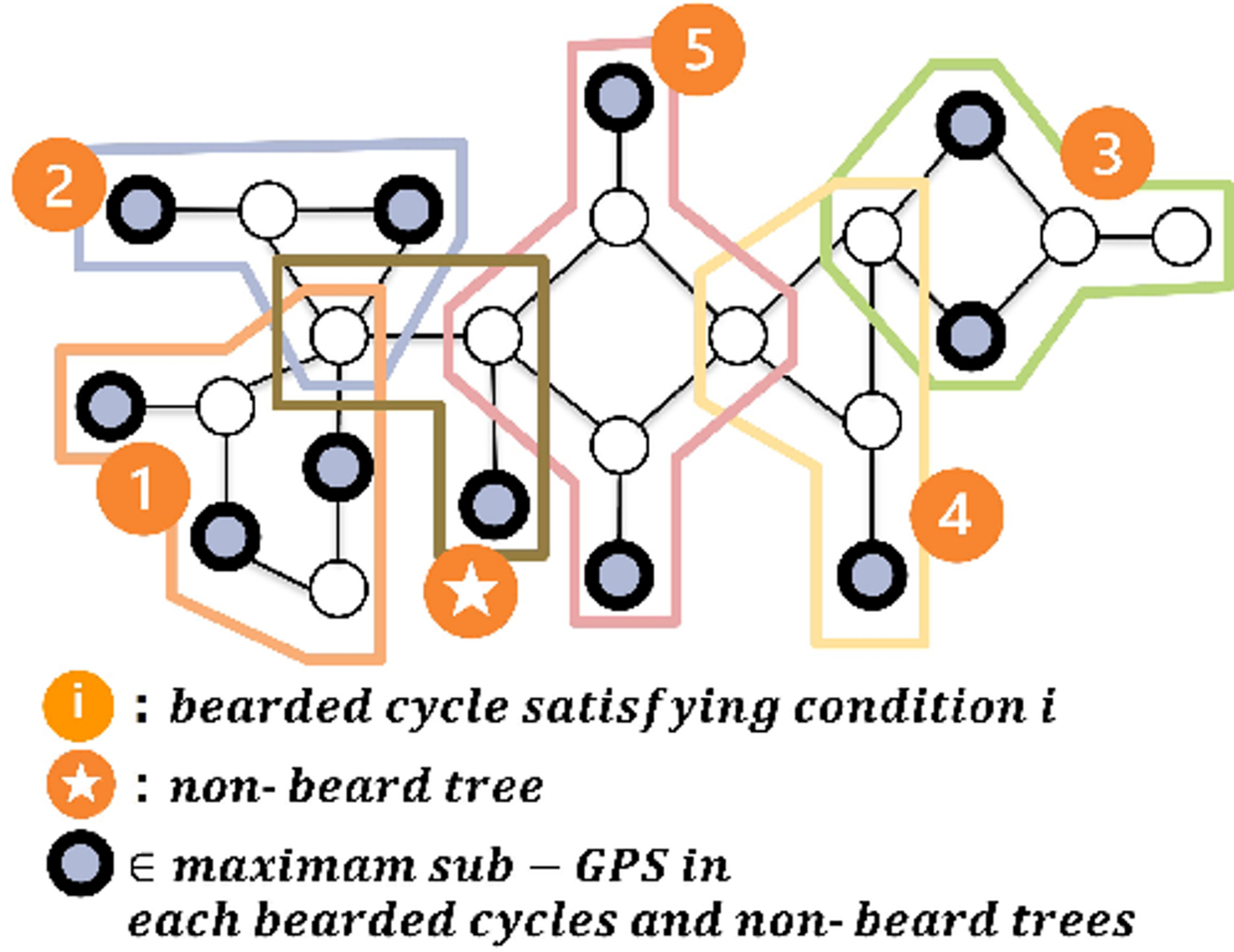}
    \caption{An example of maximum $\GPS$ consisting of $\subGPS$s in subgraphs}
    \label{fig:exp_subgpses}
    \end{center}
\end{figure}

In addition, the size of the maximum $\subGPS$ for a \nt~is equal to $L_{G_S}$. Figures \ref{fig:beardcycle_subGPS_condition1} to \ref{fig:unbeardtree_subGPS} illustrate examples of maximum $\subGPS$s for each type of \bc~and \nt.
An example of the maximum $\subGPS$ in a cactus graph is depicted in Figure \ref{fig:exp_subgpses}.

Finally, we prove the union of $\subGPS$s of all \bc s and \nt s in a non-cycle cactus graph $G$ is a maximum $\GPS$ in $G$.

\begin{theorem*}[\textbf{Theorem \ref{thm:main}}]
Let $G$ be a cactus graph that is not a simple cycle.\\
Let $BC_{all} = \{ BC_0, BC_1, \ldots, BC_k \}$ be the set of all \bc s in $G$. 
Then, the general position number $\GPN(G)$ is:
\begin{equation}
\notag
    \GPN(G) = 2n_1 + n_2 + n_3 + |L_G|
\end{equation}
where,
\begin{itemize}
  \item $n_1=|\{ BC\in BC_{all}\mid |P_{NA}(BC)|>|P_{A}(BC)|\}|$
  \item $n_2=|\{ BC\in BC_{all}\mid 0<|P_{NA}(BC)|\le |P_{A}(BC)|\}|$
  \item $n_3=|\{ BC\in BC_{all}\mid |P_{NA}(BC)|=0 \wedge |P'_{NA}(BC)|>|P'_{A}(BC)|\}|$
\end{itemize}
\end{theorem*}
\begin{proof}
Let $S$ be the union of maximum $\subGPS$s of all \bc s and \nt s os $G$. By Lemma \ref{lem:unionsubGPS}, $S$ is a $\GPS$ of $G$.
Assume, for contradiction, that there exists a $\GPS$ $S'$ of $G$ with $|S'|>|S|\wedge S'\cap A_G =\emptyset$.  We divide the discussion based on whether $G$ is a \bc~ or a \nt .

\begin{enumerate}[\arabic*.]

\item
Case $G$ is a \bc~ or a \nt~ : 
By Lemmas \ref{lem:cond1},...,\ref{lem:nonbt}, $S$ is a maximum $\GPS$ in $G$, therefore no $\GPS$ $S'$ with $|S'|>|S|$ exists, a contradiction.

\item
Case $G$ is neither a \bc~ nor a \nt~ : 
Since $G$ contains at least two \bc s or \nt s and by Lemma \ref{lem:subGPSAroundGS}, we have $|S|\ge 4$.
From $|S'|>|S|$, there must exist some \bc~or \nt~$G_S=(V_S,E_S)$ in $G$ such that $|V_S\cap S'|>|V_S\cap S|$. 
Let $G_{NS}=G-V_{S}\setminus (A_{G}\setminus B_{G_{S}})-E_{S}$, which will be used in the sequel. 
We divide the discussion based on the cases where $G_S$ is a \nt~or a \bc.

\begin{enumerate}[\alph*.]

\item
Case $G_S$ is a \nt~: 
Because $S'\cap A_G=\emptyset$ and $|V_S\cap S|=|L_{G_S}|$ and $V_S\setminus L_G\subseteq A_G$, we have $(V_S\setminus L_G)\cap S'=\emptyset$. Hence $|V_S\cap S'|\le |L_{G_S}|$, contradicting $|V_S\cap S'|>|V_S\cap S|=|L_{G_S}|$.

\item
Case $G_S$ is a \bc~: 
First, $S'\cap A_G=\emptyset$ implies $(B_{G_S}\setminus C_{G_S})\cap S'\subseteq L_{G_S}$. We divide the discussion based on the value of $|V_S\cap S'|$.

\begin{enumerate}[\roman*.]

\item
Case $|V_S\cap S'|\ge 4+|L_{G_S}|$ : 
If $|C_{G_S}\cap S'|\ge 4$, then $C_{G_S}\cap S'$ contradict the definition of $\GPS$. 
If $|(B_{G_S}\setminus C_{G_S})\cap S'|\ge |L_{G_S}|$, then, contradicting $S'\cap A_G=\emptyset$, or, for any $b\in (B_{G_S}\setminus C_{G_S})\cap S'$ there exists $c\in C_{G_S}\cap S'$ such that some shortest path from $b$ to $c$ contains at least one other vertex of $C_{G_S}\cap S'$, contradicting the definition of $\GPS$. 

\item
Case $|V_S\cap S'|= 3+|L_{G_S}|$ : 
If $|L_{G_S}|\ge 1$, then contradicting by the same argument in 2.b.i.. 
If $|L_{G_S}|=0$, then for any vertex $n$ in $G-V_S$ there exists some $s\in V_S\cap S'$ such that some shortest path from $n$ to $s$ contains another vertex in $V_S\cap S'$. Hence $G-V_S$ contains no vertex of $S'$, so $|S'|=3$, contradicting $|S'|>|S|\ge 4$.

\item
Case $|V_S\cap S'|= 2+|L_{G_S}|$ : 
We divide the discussion based on $P_{NA}(G_S)$ and $P_{A}(G_S)$.

\begin{enumerate}[\Alph*.]
    
\item
Case $|P_{NA}(G_S)|>|P_{A}(G_S)|$ : 
From lemma \ref{lem:cond1}, $|V_S\cap S|= 2+|L_{G_S}|$, this contradicts $|V_S\cap S'|>|V_S\cap S|$

\item
Case $|P_{NA}(G_S)|\le|P_{A}(G_S)|$ : 
From $|V_S\cap S'|= 2+|L_{G_S}|$ and $(B_{G_S}\setminus C_{G_S})\cap S'\subseteq L_{G_S}$, we deduce $|C_{G_S}\cap S'|\ge 2$.
We divide the discussion based on $|C_{G_S}\cap S'|$.

\begin{description}

\item[I.]
Case $|C_{G_S}\cap S'|\ge 3$ : 
For any $b\in (B_{G_S}\setminus C_{G_S})\cap S'$ there exists $c\in C_{G_S}\cap S'$ such that some shortest path from $b$ to $c$ contains at least one other vertex of $C_{G_S}\cap S'$. Hence $(B_{G_S}\setminus C_{G_S})\cap S'=\emptyset$, therefore $|S'|=|C_{G_S}\cap S'|=3$(note that $|C_{G_S}\cap S'|$ contradict the definition of $\GPS$), contradicting $|S'|>|S|\ge 4$.

\item[II.]
Case $|C_{G_S}\cap S'|=2$ : 
Let $S'_{NC}=(V\setminus C_{G_S})\cap S'$ and $P_{S'_{NC}}=V_{mP}(G_S,S'_{NC})\cap C_{G_S}$. By the definition of $\GPS$ and $C_{G_S}\cap S'\neq\emptyset$, we have $C_{G_S}\setminus P_{S'_{NC}}\neq \emptyset$ and $C_{G_S}\cap S'\subseteq C_{G_S}\setminus P_{S'_{NC}}$.
Consider four vertices $p_{1},p_{2}\in P_{S'_{NC}}$ and $n_{1},n_{2}\in C_{G_{S}}\setminus P_{S'_{NC}}$ with $(p_{1},n_{1}),(p_{2},n_{2})\in E$. 
Furthermore, consider two vertices $s_{1},s_{2}\in C_{G_{S}}\cap S'$ such that no shortest path from $s_1$ to $n_1$ contains $s_2$ and no shortest path from $s_2$ to $n_2$ contains $s_1$ (these exist by the definition of $\GPS$). 
Let $d_{s_{1}n_{1}}$, $d_{s_{1}n_{2}}$, $d_{n_{1}n_{2}}$, $d_{p_{1}p_{2}}$ denote : 
\begin{description}
    \item[$d_{s_{1}n_{1}}$] : length of the simple path from $s_1$ to $n_1$ in $G_S$ not containing $s_2$
    \item[$d_{s_{1}n_{2}}$] : length of the simple path from $s_1$ to $n_2$ in $G_S$ containing $s_2$
    \item[$d_{n_{1}n_{2}}$] : length of the shortest path from $n_1$ to $n_2$ in $G_S$
    \item[$d_{p_{1}p_{2}}$] : length of the simple path from $p_1$ to $p_2$ in $G_S$ not containing $s_1,s_2$
\end{description}
If $d_{n_{1}n_{2}}\ge d_{p_{1}p_{2}}+2$ then $|C_{G_{S}}\setminus P_{S'_{NC}}|>|P_{S'_{NC}}|$. 
If $d_{n_{1}n_{2}}<d_{p_{1}p_{2}}+2$ then $d_{n_{1}n_{2}}=d_{n_{1}s_{1}}+d_{s_{1}n_{2}}$.
From the definition of $\GPS$, $d_{s_{1}n_{2}}>d_{s_{1}n_{1}}+d_{p_{1}p_{2}} $.
Hence $d_{n_{1}n_{2}}>2d_{s_{1}n_{1}}+d_{p_{1}p_{2}}$, whitch implies $d_{n_{1}n_{2}}>d_{p_{1}p_{2}}$, so $|C_{G_{S}}\setminus P_{S'_{NC}}|>|P_{S'_{NC}}|$.
Next, we show that there exists a connected component of $G_{NS}$ that contains no vertex of $S'$.
From the definitions of $P_{S'_{NC}}$ and $P_{A}(G_{S})$, we have $P_{S'_{NC}}\subseteq P_{A}(G_{S})$.
From $|P_{NA}(G_{S})|\le |P_{A}(G_{S})|\wedge P_{S'_{NC}}\subseteq P_{A}(G_{S})\wedge |C_{G_{S}}\setminus P_{S'_{NC}}|>|P_{S'_{NC}}|$, we have $|(C_{G_{S}}\setminus P_{S'_{NC}})\cap A_{G}|\ge 1$. 
Because $|V_S\cap S'|=2+|L_{G_S}|$ and $|C_{G_S}\cap S'|=2$, we have $|(V_S\setminus C_{G_S})\cap S'|=|L_{G_S}|$, so $(C_{G_S}\setminus P_{S'_{NC}})\cap A_G\subseteq (C_{G_S}\cap A_G)\setminus B_{G_S}$. By this inclusion and the definition of $P_{S'_{NC}}$, it follows that there exists a connected component of $G_{NS}$ that contains no vertex of $S'$.
Let $G_E=(V_E,E_E)$ be a connected component of $G_{NS}$ that contains no vertex of $S'$. We now show $(|V_{S}\cap S|+|V_{E}\cap S|)-(|V_{S}\cap S'|+|V_{E}\cap S'|)\ge 0$.
By lemma \ref{lem:subGPSAroundGS} we have $|V_E\cap S|\ge 2$. From $|P_{NA}(G_S)|\le |P_{A}(G_S)|$ and Lemmas \ref{lem:cond1},...,\ref{lem:nonbt} it follows that $|L_{G_S}|\le |V_E\cap S|\le 1+|L_{G_S}|$. 
Using $|V_E\cap S|\ge 2$, $|V_E\cap S|=0$, $|V_S\cap S'|= 2+|L_{G_S}|$, $|L_{G_S}|\le |V_E\cap S|$, we get $(|V_{S}\cap S|+|V_{E}\cap S|)-(|V_{S}\cap S'|+|V_{E}\cap S'|)\ge (2+|L_{G_S}|)-(2+|L_{G_S}|)=0$.
From this, contradicting $|S'|>|S|$.

\end{description}
\end{enumerate}
\item
Case $|V_S\cap S'|= 1+|L_{G_S}|$ : 
We divide the discussion based on the sizes of $P_{A}(G_{S})$,$P_{NA}(G_{S})$,$P'_{A}(G_{S})$,$P'_{NA}(G_{S})$.

\begin{enumerate}[\arabic*.]

\item
Case $|P_{NA}(G_{S})|>0$ or $|P_{NA}(G_{S})|=0\wedge |P'_{NA}(G_{S})|>|P'_{A}(G_{S})|$ : 
By lemmas \ref{lem:cond1},\ref{lem:cond2},\ref{lem:cond3}, we have $|V_S\cap S|\ge 1+|L_{G_S}|$, contradicting $|V_S\cap S'|>|V_S\cap S|$.

\item
Case $|P_{NA}(G_{S})|=0\wedge 0<|P'_{NA}(G_{S})|\le |P'_{A}(G_{S})|$ or  $|P'_{NA}(G_{S})|=0$ : 

By lemmas \ref{lem:cond4},\ref{lem:cond5}, we have $|V_S\cap S| = |L_{G_S}|$.
Since $(V_{S}\setminus C_{G_{S}})\cap S'\subseteq L_{G_{S}}$, it follows that $|(V_{S}\setminus C_{G_{S}})\cap S'|\le |L_{G_{S}}|$. We divide the discussion based on the value of $|C_{G_{S}}\cap S'|$.

\begin{description}

\item[I.]
Case $| C_{G_{S}}\cap S'|\ge 3$ : 
Then for any $b\in  (B_{G_{S}}\setminus C_{G_{S}})\cap S'$ there exists $c\in C_{G_{S}}\cap S'$ such that some shortest path from $b$ to $c$ contains another vertex of $C_{G_{S}}\cap S'$. Hence $(B_{G_{S}}\setminus C_{G_{S}})\cap S'=\emptyset$ and $|S'|=|C_{G_{S}}\cap S'|=3$, contradicting $|S'|>|S|\ge 4$.

\item[II.]
Case $|C_{G_{S}}\cap S'|=2$ : 
Then $|(V_{S}\setminus C_{G_{S}})\cap S'|=|L_{G_{S}}|-1$. 
Let $(V\setminus C_{G_{S}})\cap S'=S'_{NC}$ and $V_{mP}(G_{S},S'_{NC})\cap C_{G_{S}}=P_{S'_{NC}}$.
By the definition of $\GPS$, we have $C_{G_{S}}\setminus P_{S'_{NC}}\neq \emptyset \wedge C_{G_{S}}\cap S'\subseteq C_{G_{S}}\setminus P_{S'_{NC}}$.
Consider four vertices $p_{1},p_{2}\in P_{S'_{NC}}$ and $n_{1},n_{2}\in C_{G_{S}}\setminus P_{S'_{NC}}$ with $(p_{1},n_{1}),(p_{2},n_{2})\in E$. 
Furthermore, consider two vertices $s_{1},s_{2}\in C_{G_{S}}\cap S'$ such that no shortest path from $s_1$ to $n_1$ contains $s_2$ and no shortest path from $s_2$ to $n_2$ contains $s_1$ (these exist by the definition of $\GPS$). 
Let $d_{s_{1}n_{1}}$, $d_{s_{1}n_{2}}$, $d_{n_{1}n_{2}}$, $d_{p_{1}p_{2}}$ denote : 
\begin{description}
    \item[$d_{s_{1}n_{1}}$] : length of the simple path from $s_1$ to $n_1$ in $G_S$ not containing $s_2$
    \item[$d_{s_{1}n_{2}}$] : length of the simple path from $s_1$ to $n_2$ in $G_S$ containing $s_2$
    \item[$d_{n_{1}n_{2}}$] : length of the shortest path from $n_1$ to $n_2$ in $G_S$
    \item[$d_{p_{1}p_{2}}$] : length of the simple path from $p_1$ to $p_2$ in $G_S$ not containing $s_1,s_2$
\end{description}
By the same argument in 2.b.iii.B.II. , we have $|(C_{G_{S}}\setminus P_{S'_{NC}})\cap A_{G}|\ge 1$.
By $|P_{NA}(G_{S})|=0\wedge |P'_{NA}(G_{S})|\le |P'_{A}(G_{S})|$ or $|P_{NA}(G_{S})|=0\wedge |P'_{NA}(G_{S})|=0$, and there exists exactly one beard containing no vertex in $S'$(by $|(V_{S}\setminus C_{G_{S}})\cap S'|=|L_{G_{S}}|-1$), it follows that there exists a connected component of $G_{NS}$ that contains no vertex in $S'$.
Let $G_E=(V_E,E_E)$ be a connected component of $G_{NS}$ that contains no vertex in $S'$. We now show $(|V_{S}\cap S|+|V_{E}\cap S|)-(|V_{S}\cap S'|+|V_{E}\cap S'|)\ge 0$.
By lemma \ref{lem:subGPSAroundGS} we have $|V_E\cap S|\ge 2$. 
From $|P_{NA}(G_{S})|=0\wedge |P'_{NA}(G_{S})|\le |P'_{A}(G_{S})|$ or $|P_{NA}(G_{S})|=0\wedge |P'_{NA}(G_{S})|=0$, it follows that $|V_{S}\cap S|=|L_{G_{S}}|$.
Using $|V_{E}\cap S|\ge 2\wedge |V_{E}\cap S'|=0\wedge |V_{S}\cap S|=|L_{G_{S}}|\wedge |V_{S}\cap S'|=1+|L_{G_{S}}|$, we get $(|V_{S}\cap S|+|V_{E}\cap S|)-(|V_{S}\cap S'|+|V_{E}\cap S'|)\ge (2+|L_{G_{S}}|)-(1+|L_{G_{S}}|)=1$.
From this, contradicting $|S'|>|S|$.

\item[III.]
Case $|C_{G_{S}}\cap S'|=1$ : 
Then $|(V_{S}\setminus C_{G_{S}})\cap S'|=|L_{G_{S}}|$.
Let $S'_{NC}=(V\setminus C_{G_S})\cap S'$ and $P_{S'_{NC}}=V_{mP}(G_S,S'_{NC})\cap C_{G_S}$. By the definition of $\GPS$ and $C_{G_S}\cap S'\neq\emptyset$, we have $C_{G_S}\setminus P_{S'_{NC}}\neq \emptyset$ and $C_{G_S}\cap S'\subseteq C_{G_S}\setminus P_{S'_{NC}}$.
From each beard of $G_S$ contains exactly one vertex in $S'$, we get $(C_{G_{S}}\setminus  P_{S'_{NC}})\cap A_{G}\subseteq (C_{G_{S}}\cap A_{G})\setminus B_{G_{S}}$, it follows that there exists a connected component of $G_{NS}$ that contains no vertex in $S'$. 
Let $G_E=(V_E,E_E)$ be a connected component of $G_{NS}$ that contains no vertex in $S'$. We now show $(|V_{S}\cap S|+|V_{E}\cap S|)-(|V_{S}\cap S'|+|V_{E}\cap S'|)\ge 0$.
By the same argument in 2.b.iv.B.II. , we have $(|V_{S}\cap S|+|V_{E}\cap S|)-(|V_{S}\cap S'|+|V_{E}\cap S'|)\ge (2+|L_{G_{S}}|)-(1+|L_{G_{S}}|)=1$.
From this, contradicting $|S'|>|S|$.

\end{description}
\end{enumerate}

\item
Case $|V_S\cap S'|=|L_{G_S}|$ : 
By lemma \ref{lem:cond1},...,\ref{lem:nonbt}, we have $|V_S\cap S|\ge |L_{G_S}|$, contradicting $|V_S\cap S'|>|V_S\cap S|$.

\end{enumerate}
\end{enumerate}
\end{enumerate}

From above, for any \bc~or \nt~in $G$, we derive a contradiction, therefore there exist no $\GPS$ $S'$ such that $|S'|>|S|$ in $G$. $S$ is a maximum $\GPS$ in $G$.

Finally, by Lemma \ref{lem:sizeofunion}, $\GPN(G)=|S|=2n_{1}+n_{2}+n_{3}+|L_{G}|$.
\end{proof}

\section{Conclusion}
In this paper, we presented the (exact) General Position Number ($\GPN$) for cactus graphs, which is the maximum size of General Position Set ($\GPS$).
Previous studies have provided exact $\GPN$, upper bounds, or lower bounds of the $\GPN$ for certain specific classes of graphs. 
For cactus graphs, both upper and lower bounds for arbitrary cactus graphs were known, 
and exact $\GPN$ had been investigated only for cactus graphs satisfying specific conditions.

In this work, we provide the exact $\GPN$ for general cactus graphs without any restrictions. 
To achieve this, we present a novel graph decomposition technique, which decomposes a given cactus graph $G$ 
into components called $\bc$s and $\nt$s which we newly defined. 
We then show that by finding $\subGPS$, which is a $\GPS$ satisfying certain conditions, in each component, 
and the union of these $\subGPS$s forms a $\GPS$ for the entire graph $G$ and furthermore, that this union is of maximum size.

As future directions, we consider extending the study to determine the $\GPS$ for other classes of graphs beyond cactus graphs, 
as well as investigating the exact value of the mutual visibility number for cactus graphs.
The techniques developed in this work are expected to be applicable to the study of $\GPN$ in other graph classes. 
Moreover, since mutual visibility sets share many similar properties with $\GPS$, 
the findings of this study may also be useful in the analysis of mutual visibility sets in cactus graphs.

\bibliographystyle{unsrt}  
\bibliography{references}

\end{document}